\documentclass[12pt,fullpage]{article}

\usepackage{amsfonts}
\usepackage{xr}
\usepackage{graphicx}
\usepackage{amsmath}
\usepackage{epsfig}
\usepackage{amssymb}

\def\Z{\mathbb{Z}}
\def\R{\mathbb{R}}
\def\N{\mathbb{N}}

\def\C{\mathbb{C}}
\def\epsilon{\varepsilon}

\def\tilde{\widetilde}

\def\trait (#1) (#2) (#3){\vrule width #1pt height #2pt depth #3pt}
\def\fin{\hfill\trait (0.1) (5) (0) \trait (5) (0.1) (0) \kern-5pt 
\trait (5) (5) (-4.9) \trait (0.1) (5) (0)}

\newcommand{\be}{\begin{equation}}
\newcommand{\ee}{\end{equation}}
\newcommand{\baa}{\begin{array}}
\newcommand{\eaa}{\end{array}}
\newcommand{\ba}{\begin{eqnarray}}
\newcommand{\ea}{\end{eqnarray}}

\newtheorem{theo}{\bf Theorem}[section]
\newtheorem{lem}[theo]{\bf Lemma}
\newtheorem{prop}[theo]{\bf Proposition}

\newtheorem{defi}[theo]{\bf Definition}
\newtheorem{rem}[theo]{\bf Remark}

\begin{document}
\date{}
\title{\bf{Sobolev and Hardy-Sobolev spaces on graphs}}
\author{Emmanuel Russ$^{\hbox{\small{ a}}}$ and Maamoun Turkawi$^{\hbox{\small{ b}}}$\\
\\
\footnotesize{$^{\hbox{a }}$Universit\'e Joseph Fourier, Institut Fourier, 100 rue des Maths, BP 74, F-38402 St-Martin d'H\`eres, France}\\
\footnotesize{$^{\hbox{b }}$Aix-Marseille Universit\'e, LATP, 
Facult\'e des Sciences et Techniques, Case cour A}\\
\footnotesize{Avenue Escadrille Normandie-Niemen, F-13397 Marseille 
Cedex 20, France}\\
}
\maketitle

\begin{abstract}
Let $\Gamma$ be a graph. Under suitable geometric assumptions on $\Gamma$, we give several equivalent characterizations of Sobolev and Hardy-Sobolev spaces on $\Gamma$, in terms of maximal functionals, Haj\l asz type functionals or atomic decompositions. As an application, we study the boundedness of Riesz transforms on Hardy spaces on $\Gamma$. This gives the discrete counterpart of the corresponding results on Riemannian manifolds.
\end{abstract}

\tableofcontents


\section{Introduction}
\subsection{The Euclidean case}
Let $n\in \N^{\ast}$ and $1\leq p\leq +\infty$. Throughout the paper, if $A(f)$ and $B(f)$ are two quantities depending on a function $f$ ranging in a set $E$, say that $A(f)\lesssim B(f)$ if and only if there exists $C>0$ such that, for all $f\in E$,
$$
A(f)\leq CB(f),
$$
and that $A(f)\sim B(f)$ if and only if $A(f)\lesssim B(f)$ and $B(f)\lesssim A(f)$. \par
\noindent The classical $W^{1,p}(\R^n)$ space, or its homogenous version $\dot{W}^{1,p}(\R^n)$, can be characterized in terms of maximal functions. Namely, if $f\in L^1_{loc}(\R^n)$, define, for all $x\in \R^n$, 
$$
Nf(x):=\sup_{B\ni x} \frac 1{r(B)\left\vert B\right\vert} \int_B \left\vert f(y)-f_B\right\vert dy,
$$
where the supremum is taken over all balls $B$ containing $x$ and 
$$
f_B:=\frac 1{\left\vert B\right\vert} \int_B f(y)dy
$$
is the mean value of $f$ over $B$. Here and after in this section, if $B\subset \R^n$ is a ball, $\left\vert B\right\vert$ stands for the Lebesgue measure of $B$ and $r(B)$ for its radius. \par
\noindent Then (\cite{calderon}), for $1<p\leq +\infty$, $\nabla f\in L^p(\R^n)$ if and only if $Nf\in L^p(\R^n)$, and
$$
\left\Vert \nabla f\right\Vert_{L^p(\R^n)}\sim \left\Vert Nf\right\Vert_{L^p(\R^n)}.
$$
Another maximal function characterizing Sobolev spaces was introduced in \cite{art}. For $f\in L^1_{loc}(\R^n)$ and $x\in \R^n$, define
$$
Mf(x):=\sup \left\vert \int_{\R^n} f(y)\mbox{div }\Phi(y)dy\right\vert,
$$
where the supremum is taken over all vector fields $\Phi\in 
L^{\infty}(\R^n,\C^n)$, whose distributional divergence is a
bounded function in
$\R^n$,
supported in a ball $B\subset \R^n$ 
containing $x$, with
\[
\left\Vert {\Phi}\right\Vert_{\infty} + r(B) \left\Vert \mbox{div }
\Phi\right\Vert_{\infty}\leq \frac 1{\left\vert B\right\vert}.
\]
Then (\cite{art}), for $1<p\leq +\infty$, $\nabla f\in L^p(\R^n)$ if and only if $Nf\in L^p(\R^n)$, and
$$
\left\Vert \nabla f\right\Vert_{L^p(\R^n)}\sim \left\Vert Nf\right\Vert_{L^p(\R^n)}.
$$
Another description of Sobolev spaces is due to Haj\l asz. For $f\in L^1_{loc}(\R^n)$, $1\leq p\leq +\infty$, say that $f\in \dot{M}^{1,p}(\R^n)$ if and only if there exists $g\in L^p(\R^n)$ such that, for all $x,y\in \R^n$,
\begin{equation} \label{fg}
\left\vert f(x)-f(y)\right\vert\leq d(x,y)(g(x)+g(y)).
\end{equation}
Set
$$
\left\Vert f\right\Vert_{\dot{M}^{1,p}(\R^n)}:=\inf \left\Vert g\right\Vert_{L^p(\R^n)},
$$
the infimum being taken over all functions $g$ such that \eqref{fg} holds. It was proved by Haj\l asz (\cite{haj}) that, for $1<p\leq +\infty$, $f\in \dot{M}^{1,p}(\R^n)$ if and only if $\nabla f\in L^p(\R^n)$ and
\begin{equation}Ê\label{equivrn}
\left\Vert f\right\Vert_{\dot{M}^{1,p}(\R^n)}\sim \left\Vert \nabla f\right\Vert_{L^p(\R^n)}.
\end{equation}
What happens in these results when $p=1$ ? The previous results break down when $p=1$, but correct substitutes involving Hardy-Sobolev spaces can be given. More precisely (see below in the introduction), $\dot{M}^{1,1}(\R^n)$ coincides with the space of locally integrable functions with gradient in the $H^1(\R^n)$ Hardy space. \par
\noindent The $H^1(\R^n)$ Hardy space is well-known to be the right substitute for $L^1(\R^n)$ for many questions in harmonic analysis. Let us recall one possible definition of $H^1(\R^n)$. Fix a function
$\varphi\in {\cal S}(\R^n)$ such that
$\int_{\R^n} \varphi(x)dx=1$. For all $t>0$, define $\varphi_t(x):=t^{-n}
\varphi\left(\frac xt\right)$. Define then $H^1(\R^n)$ as the space of locally integrable functions $f$ on $\R^n$ such that the vertical maximal function
\[
{\mathcal M}f(x):=\sup\limits_{t>0} \left\vert \varphi_t \ast f(x)\right\vert
\]
belongs to $L^1(\R^n)$. Define
\[
\left\Vert f\right\Vert_{H^1(\R^n)} := \left\Vert {\mathcal
M}f\right\Vert_{L^1(\R^n)}.
\]
As for classical Sobolev spaces, let us consider the Hardy-Sobolev space $H^{1,1}(\R^n)$ made of functions $f\in L^1(\R^n)$ such that $\nabla f\in H^1(\R^n)$, in the sense that, for all $1\leq j\leq n$, $\frac{\partial f}{\partial x_j}\in H^1(\R^n)$. Define also $\dot{H}^{1,1}(\R^n)$ as the space of functions $f\in L^1_{loc}(\R^n)$ such that $\nabla f\in H^1(\R^n)$, equipped with the semi-norm
$$
\left\Vert f\right\Vert_{\dot{H}^{1,1}(\R^n)}:=\left\Vert \nabla f\right\Vert_{H^1(\R^n)}.
$$ 
\noindent Various characterizations of this space (as well as its adaptations to the case of domains of $\R^n$) were given in the literature. It can be described in terms of a functional involving second order differences (\cite{stri}). In \cite{miyachi}, $H^{1,1}(\R^n)$ was characterized in terms of the maximal function $Nf$. Namely, for $f\in L^1_{loc}(\R^n)$, $\nabla f\in H^1(\R^n)$ if and only if $Nf\in L^1(\R^n)$ and
$$
\left\Vert Nf\right\Vert_{L^1(\R^n)}\sim \left\Vert \nabla f\right\Vert_{H^1(\R^n)}:=\sum_{j=1}^n \left\Vert \frac{\partial f}{\partial x_j}\right\Vert_{H^1(\R^n)}.
$$
It was shown in \cite{art} that the functional $Mf$ defined above characterizes Hardy-Sobolev spaces (actually, this was the reason why this maximal function was introduced in \cite{art}, since it is particularly suited to the study of Hardy-Sobolev spaces on strongly Lipschitz domains of $\R^n$). More precisely, $\nabla f\in H^1(\R^n)$ if and only if $Mf\in L^1(\R^n)$ and
$$
\left\Vert Mf\right\Vert_{L^1(\R^n)}\sim \left\Vert \nabla f\right\Vert_{H^1(\R^n)}.
$$
Moreover, going back to Haj\l asz's functional, it was proved in \cite{kosksaks} that $f\in \dot{M}^{1,1}(\R^n)$ if and only if $\nabla f\in H^1(\R^n)$ and
$$
\left\Vert f\right\Vert_{\dot{M}^{1,1}(\R^n)}\sim \left\Vert \nabla f\right\Vert_{H^1(\R^n)}.
$$
Finally, an atomic decomposition for Hardy-Sobolev spaces was given in \cite{stri}. In this paper, an atom is a function $b$ supported in a cube such that $(-\Delta)^{1/2}b$ satisfies suitable $L^p$ estimates (\cite{stri}, definition 5.1).\par
\noindent Another characterization of $H^1(\R^n)$ states that it is exactly the space of functions $f\in L^1(\R^n)$ such that, for all $1\leq j\leq n$, $\frac{\partial}{\partial x_j}(-\Delta)^{-1/2}f\in L^1(\R^n)$ (see \cite{feffstein}). The operators $R_j:=\frac{\partial}{\partial x_j}(-\Delta)^{-1/2}f$ are the Riesz transforms. Thus, $(-\Delta)^{-1/2}$ maps continuously $H^1(\R^n)$ into $\dot{H}^{1,1}(\R^n)$. 
\subsection{The case of Riemannian manifolds}
These various characterizations can be extended to the framework of Riemannian manifolds. Namely, let $M$ be a complete Riemannian manifold, endowed with its Riemannian metric $d$ and its Riemannian measure $\mu$. Say that $M$ satisfies the doubling condition if there exists $C>0$ such that, for all $x\in M$ and all $r>0$,
$$
\mu(B(x,2r))\leq C\mu(B(x,r)).
$$
Say that $M$ satisfies an $L^1$ scaled Poincar\'e inequality on balls if there exists $C>0$ such that, for all balls $B\subset M$ with radius $r$ and all functions $f\in C^{\infty}(B)$,
$$
\int_B \left\vert f(x)-f_B\right\vert d\mu(x)\leq Cr\int_B \left\vert df(x)\right\vert d\mu(x).
$$
Define the $\dot{M}^{1,p}$ spaces and the $Nf$ functional as in the Euclidean case. Then, for $1\leq p<+\infty$, $f\in \dot{M}^{1,p}$ if and only if $Nf\in L^p(M)$ (\cite{kin2007}). A version of the maximal function in \cite{art} is given in \cite{maxhs}, where it is shown that it characterizes $\dot{M}^{1,1}$. Moreover, an atomic decomposition for $\dot{M}^{1,1}$ is provided in \cite{badr2010atomic}, where it is also shown that $f\in \dot{M}^{1,1}$ if and only if $df$ belongs to the Hardy space of exact differential forms $H^1_{d}(\Lambda^1T^{\ast}M)$ introduced in \cite{amr}. Since $d\Delta^{-1/2}$ is bounded from $H^1_{d^{\ast}}(\Lambda^0T^{\ast}M)$ from $H^1_{d}(\Lambda^1T^{\ast}M)$ (see \cite{amr}, Theorem 5.16), if $\Delta$ denotes the Laplace-Beltrami operator, $\Delta^{-1/2}$ maps continuously $H^1_{d^{\ast}}(\Lambda^0T^{\ast}M)$ into $\dot{M}^{1,1}$. \par

\medskip

\noindent In the present work, we investigate Sobolev and Hardy-Sobolev spaces on graphs, and establish the discrete counterpart of the results obtained on Riemannian manifolds. Namely, we characterize Sobolev and Hardy-Sobolev spaces in terms of maximal functions and provide an atomic decomposition for Hardy-Sobolev spaces. We also investigate the boundedness of Riesz transforms on Hardy spaces. 
\section{Description of the results} \label{pres}
\subsection{Presentation of the graph} \label{graph}
The geometric context is the same as in \cite{BR09}, and we recall it for the sake of completeness. Let $\Gamma$ be an infinite set and $\mu_{xy}=\mu_{yx}$ a symmetric weight on $\Gamma \times \Gamma$. Say that $x\sim y$ if and only if $\mu_{xy}>0$, and let $E$ stand for the set of edges in $\Gamma$, defined as the set of $(x,y)\in \Gamma\times \Gamma$ such that $\mu_{xy}>0$. For all $x\in \Gamma$, say that $x$ is a vertex of $\Gamma$. \par
\noindent For $x,y\in \Gamma$, a path joining $x$ to $y$ is a finite sequence of vertices $x_{0}=x,\cdots,x_{N}=y$ such that, for all $0\leq i \leq N-1, x_{i}\sim x_{i+1}$.
Say that this path has length $N$. Assume that $\Gamma$ is connected, which means that, for all $x,y\in \Gamma,$ there exists a path joining $x$ to $y$. The distance between $x$ and $y$, denoted $d(x,y)$, is defined as the shortest length of a path joining $x$ and $y$. For all $x\in \Gamma$ and all $r\geq 0$, define the closed ball $$B(x,r):=\{y\in\Gamma; d(x,y)\leq r\}.$$
\noindent In the sequel, we always assume that $\Gamma$ is locally uniformly finite, which means that there exists $N\in\mathbb{N}^{\ast}$ such that, for all $x\in\Gamma,\# B(x,r)\leq N$.\par
\noindent For any subset $\Omega\subset \Gamma,$
set
$$\partial \Omega:=\{x\in \Omega; \exists y\sim x, y\notin  \Omega \}$$
and
$$
\overset{\circ}{\Omega}:=\Omega\setminus \partial\Omega.
$$
In other words, $\overset{\circ}{\Omega}$ is the set of points $x\in  \Omega$ such that $y\in \Omega$ whenever $x\sim y$.
Denote by $E_{\Omega}$ the set of edges in $\Omega$,
$$E_{\Omega}=\{(x,y)\in \Omega \times \Omega: x\sim y,x,y\in \Omega\}.$$
\noindent We also define a distance on $E$. For $\gamma=(x,y)$ and $\gamma^{\prime}=(x^{\prime},y^{\prime})\in E$, set
$$
d(\gamma,\gamma^{\prime}):=\max\left(d(x,x^{\prime}),d(y,y^{\prime})\right).
$$
\subsubsection{The measures on $\Gamma$ and $E$} \label{measures}
For all $x\in \Gamma,$ set $m(x)=\sum \limits_{y\sim x} \mu_{xy}$ (recall that this sum has at most $N$ terms). We always assume in the sequel that $m(x)>0$ for all $x\in \Gamma.$   If $\Omega \subset \Gamma,$ define $m(\Omega)=\sum \limits_{x\in \Omega}m(x)$. For all $x\in \Gamma$ and $r> 0$, write $V(x,r)$ instead of $m(B(x,r))$ and, if $B$ is a ball, $m(B)$ will be denoted by $V(B)$.\par
\noindent Here is a growth assumption on the volume of balls of $\Gamma$, which may be satisfied or not.
\begin{defi}\label{doubling condition}
[Doubling property] Say that $(\Gamma,d,m)$ satisfies the  doubling property if there exists a constant $C> 0$ such that for all balls $B(x,r), x\in \Gamma, r>0,$
\begin{equation} \label{doubling} \tag{$D$}
V(x,2r)\leq C V(x,r).
\end{equation}
\end{defi}
This means that $(\Gamma,d,m)$ is a space of homogeneous type in the sense of Coifman and Weiss (\cite{coifmanweiss}). It is plain to check that, if $\Gamma$ satisfies \eqref{doubling}, then there exist $C,s>0$ such that, for all $x\in \Gamma$, all $r>0$
and all $\theta\geq1$,
\begin{equation} \label{defs}
V(x,\theta r) \leq C \theta^s V(x,r).
\end{equation}
\begin{rem} \label{mart}
Observe also that, since $\Gamma$ is infinite, it is also unbounded (since it is locally uniformly finite) so that, if \eqref{doubling} holds, then $m(\Gamma)=+\infty$ (see \cite{martell01}).
\end{rem}
For all $1\leq p<+ \infty$, say that a function $f:\Gamma \rightarrow \mathbb{R}$   belongs to $L^{p}(\Gamma)$ if
$$\|f\|_{_{L^{p}(\Gamma)}}=\left(\sum_{x\in \Gamma}|f(x)|^{p}m(x)\right)^{1/p}<+\infty.$$
Note that the $L^2(\Gamma)$-norm derives from the scalar product
$$
\langle f,g\rangle_{L^2(\Gamma)}:=\sum_{x\in \Gamma} f(x)g(x)m(x).
$$
Say that $f\in L^{\infty}(\Gamma)$ if
$$\|f\|_{_{L^{\infty}(\Gamma)}}=\sup_{x\in\Gamma}|f(x)|<+\infty.$$
If $B\subset \Gamma$ is a ball, denote by $L^p_0(B)$ the subspace of $L^p(\Gamma)$ made of functions $f$ supported in $B$ and satisfying 
$$
\sum_{x\in B} f(x)m(x)=0.
$$
We also need a measure on $E$. For any subset $A\subset E$, define
$$
\mu(A):=\sum_{(x,y)\in A} \mu_{xy}.
$$
It is easily checked (\cite{BR09}, Section 8) that, if \eqref{doubling} holds, then $E$, equipped with the distance $d$ and the measure $\mu$, is a space of homogeneous type.\par
\noindent Define $L^{p}$ spaces on $E$ in the following way. For $1\leq p<+\infty$, say that a function $F$ on $E$ belongs to $L^{p}(E)$ if and only if $F$ is antisymmetric, which means that $F(x,y)=-F(y,x)$ for all $(x,y)\in E$, and
$$\|F\|^{p}_{L^{p}(E)}:=\frac{1}{2}\sum_{(x,y)\in E}|F(x,y)|^{p}\mu_{xy}<+ \infty.$$
Observe that the $L^{2}(E)-$norm derives from the scalar product
$$\langle F,G\rangle_{L^{2}(E)}:=\frac{1}{2}\sum_{x,y\in \Gamma}F(x,y)G(x,y)\mu_{xy}.$$
Finally, say that $F\in L^{\infty}(E)$ if and only if $F$ is antisymmetric and
$$\|F\|_{L^{\infty}(E)}:=\frac{1}{2} \sup_{(x,y)\in E}|F(x,y)|<+\infty.$$
Define $L^p(E_{\Omega})$ similarly.
\subsubsection{The Markov kernel}
Define $p(x,y)=\frac{\mu_{xy}}{m(x)}$ for all $x,y\in \Gamma$. Observe that $p(x,y)=0$ 
if $d(x,y)\geq2$. Moreover, for all $x\in \Gamma$,
\begin{equation}\label{1}
\sum_{y\in \Gamma} p(x,y)=1
\end{equation}
and for all $x,y\in \Gamma$,
\begin{equation}\label{reversibility}
p(x,y) m(x)=p(y,x) m(y).
\end{equation}
Another assumption on $(\Gamma,\mu)$ which will be used in the sequel is a uniform lower bound for $p(x,y)$
when $x\sim y$. For $\alpha>0,$ say that $(\Gamma,\mu)$ 
satisfies the condition $\Delta (\alpha)$ if, for all $x,y\in\Gamma$,
\begin{equation}\label{deltaalpha} \tag{$\Delta(\alpha)$}
(x\sim y \Leftrightarrow \mu_{xy} \geq \alpha m(x)) \ \mbox{and}  \ x\sim x.
\end{equation}
For all functions $f$ on $\Gamma$ and all $x\in \Gamma$, define
$$Pf(x)=\sum_{y\in\Gamma}p(x,y)f(y).$$
It is easily checked (\cite{BR09}), using \eqref{reversibility}, that, for all functions $f$ on $\Gamma$,
\begin{equation} \label{I-P}
\langle (I-P)f,f\rangle=\displaystyle \frac{1}{2}\sum\limits_{x,y}p(x,y)|f(x)-f(y)|^{2}m(x).
\end{equation}
Identity \eqref{I-P} leads to the definition of the operator ``length of the gradient'' by
$$\nabla f(x)=\left(\frac{1}{2}\sum_{y\in \Gamma}p(x,y)|f(y)-f(x)|^{2}\right)^{1/2},$$
so that, for all functions $f$ on $\Gamma$,
\begin{equation} \label{bypartlaplacian}
\langle (I-P)f,f\rangle_{L^2(\Gamma)}=\left\Vert \nabla f\right\Vert_{L^2(\Gamma)}^2.
\end{equation}
\subsubsection{The differential and divergence operators}
We now define a discrete differential, following the definitions of \cite{BR09} but dealing with functions defined on subsets of $\Gamma$. Let $\Omega\subset \Gamma$. For any function $f:\Omega\rightarrow \R$ and any $\gamma=(x,y)\in E_{\Omega}$, 
define
\begin{equation}\label{df}
df(\gamma)=f(y)-f(x).
\end{equation}
The function $df$ is clearly antisymmetric on $E_{\Omega}$. Moreover, it is easily checked (\cite{BR09}, p.313) 
that, if $(\Delta(\alpha))$ holds, then for all $p\in \lbrack 1,+\infty\rbrack$ and all functions $f$ on $\Gamma$, 
\begin{equation}\label{la relation entre diff f et grad f}
\Vert df \Vert_{L^p(E)} \sim \Vert \nabla f \Vert_{L^p (\Gamma)}.
\end{equation}
We define now a divergence operator in such a way that a discrete integration by parts formula holds (see \cite{BR09}). Let $F$ be any (antisymmetric) function in $L^{2}(E_{\Omega})$. If $f$ is a function on $
\Omega$ vanishing on $\partial\Omega$ such that $df\in L^{2}(E_{\Omega})$, one has
\begin{equation*}
\begin{array}{lll}
\displaystyle \langle df,F\rangle_{L^{2}(E_{\Omega})} & = & \displaystyle \frac 12\sum_{x,y\in \Omega,\ x\sim y} df(x,y)F(x,y)\mu_{xy}\\
&  = & \displaystyle - \sum_{x,y\in \Omega,\ x\sim y} f(x)F(x,y)\mu_{xy}\\
& = & \displaystyle -\sum_{x\in \overset{\circ}{\Omega}} f(x)\left(\sum_{y\sim x,\ y\in \Gamma} p(x,y)F(x,y)\right)m(x),
\end{array}
\end{equation*}
where the second line is due to the fact that $F$ is antisymmetric and the third one holds because $f(x)=0$ when $x\in \partial\Omega$ and all the neighbours of $x$ in $\Gamma$ actually belong to $\Omega$ when $x\in\overset{\circ}{\Omega}$. Thus, if we define the divergence of $F$ by
$$ \delta F(x):=\sum_{y\sim x,\ y\in \Gamma} p(x,y)F(x,y) $$
for all $x\in \overset{\circ}{\Omega}$, it follows that
\begin{equation}\label{IPP}
\langle df ,F\rangle_{L^{2}(E_{\Omega})}=-\langle f ,\delta F\rangle_{L^{2}(\overset{\circ}{\Omega})}.
\end{equation}
\begin{rem}
A slightly different integration by parts formula on graphs can be found in \cite{cgz}, formula 2.4.
\end{rem}
\subsubsection{The Poincar\'e inequality on balls}
\begin{defi}\label{P1}
[$L^p$ Poincar\'e inequality on balls] Let $p\in [1,+\infty)$. Say that $\Gamma$ satisfies an $L^p$ scaled Poincar\'e inequality on balls if there exists a constant $C>0$
such that, for all functions $f$ on $\Gamma$ and all balls $B\subset \Gamma$ of radius $r>0$,
\begin{equation}\label{p1}\tag{$P_{p}$}
 \sum\limits_{x\in B}|f(x)-f_{B}|^pm(x)\leq C r^p \sum\limits_{x\in B} \left\vert \nabla f(x)\right\vert^p m(x),
\end{equation}
where
\begin{equation} \label{mean}
f_{B}=\frac{1}{V(B)}\sum\limits_{x\in B} f(x)m(x).
\end{equation}
\end{defi}
\begin{rem} \label{improvepoinc}
\begin{itemize}
\item[$1.$]
Note that, if $(P_1)$ holds, then one has an $L^p$ Poincar\'e inequality for all $p\in [1,+\infty)$ (see \cite{sobmetpoinc}). 
\item[$2.$]
Moreover, if $(P_p)$ holds for some $p\in (1,+\infty)$, there exists $q<p$ such that $(P_q)$ still holds (\cite{kz}). 
\end{itemize}
\end{rem}
\subsection{Sobolev spaces}
Let $\Gamma$ be a graph as in Section \ref{graph}. Let $1\leq p\leq +\infty$. Say that a scalar-valued function $f$ on $\Gamma$ belongs to the Sobolev space $W^{1,p}(\Gamma)$ if and only if
$$\left\Vert f \right\Vert_{W^{1,p}(\Gamma)}:=\left\Vert f\right\Vert_{L^{p}(\Gamma)}+\left\Vert \nabla f\right\Vert_{L^{p}(\Gamma)}<+\infty.$$
As in \cite{BR09} we will also consider the homogeneous versions of Sobolev spaces. Define $\dot{W}^{1,p}(\Gamma)$ as the space of all scalar-valued functions $f$ on
 $\Gamma$ such that $\nabla f\in L^p(\Gamma)$, equipped with the 
 semi-norm 
 $$\left\Vert f\right\Vert_{\dot{W}^{1,p}(\Gamma)}:= \left\Vert \nabla f \right\Vert
 _{L^p(\Gamma)}.$$
If $B$ is any ball in $\Gamma$ and $1\leq p \leq +\infty,$ denote by $W^{1,p}_{0}(B)$ the subspace of $W^{1,p}(\Gamma)$ made of functions supported in $\overset{\circ}B$.
\subsection{Characterizations of Sobolev spaces} \label{caracsobol}
In the present section, we give various characterizations of Sobolev spaces on graphs. The first one is formulated in terms of Haj\l asz's functionals (see \cite{hajlasz,sobmetpoinc}):
\begin{defi}\label{Hajlasz}
Let $1\leq p\leq +\infty$. 
\begin{itemize}
\item[$1.$]
The inhomogeneous Sobolev space $M^{1,p}(\Gamma)$ is defined as the space of all functions $f\in L^{p}(\Gamma)$ such that there exists a non-negative function $g\in L^{p}(\Gamma)$ satisfying
\begin{equation}\label{Hajlasz space}
|f(x)-f(y)|\leq d(x,y) \left(g(x)+g(y)\right) \mbox{ for all }x,y\in \Gamma.
\end{equation}
We equip $M^{1,p}(\Gamma)$ with the norm
\begin{equation} \label{m1pnorm}
|\vert f |\vert_{M^{1,p}(\Gamma)}:=\left\Vert f\right\Vert_{L^p(\Gamma)}+ \inf\limits_{g} |\vert g|\vert_{L^p(\Gamma)},
\end{equation}
where the infimum is taken over all functions $g\in L^p(\Gamma)$ such that \eqref{Hajlasz space} holds. 
\item[$2.$]
The homogeneous Sobolev space $\dot{M}^{1,p}(\Gamma)$ is defined as the space of all functions $f$ on $\Gamma$ such that there exists a non-negative function $g\in L^{p}(\Gamma)$ satisfying \eqref{Hajlasz space}. We equip $\dot{M}^{1,p}(\Gamma)$ with the semi-norm
$$|\vert f |\vert_{\dot{M}^{1,p}(\Gamma)}=\inf\limits_{g} |\vert g|\vert_{L^p(\Gamma)},$$
where the infimum is taken over all functions $g\in L^p(\Gamma)$ such that \eqref{Hajlasz space} holds. 
\end{itemize}
\end{defi}
\begin{rem} \label{m1pleq1}
If $B\subset \Gamma$ is a ball, define $M^{1,p}(B)$ and $\dot{M}^{1,p}(B)$, replacing $\Gamma$ by $B$ in Definition \ref{Hajlasz}.
\end{rem}
We will also characterize Sobolev spaces in terms of two maximal functions. \par
\noindent The first maximal function is modelled on the one in \cite{calderon}. For all functions $f$ on $\Gamma$ and all $x\in \Gamma$, define $Nf(x)$ by 
\begin{equation}\label{Nf}
Nf(x):=\sup_{B\ni x} \frac{1}{r(B)V(B)}\sum_{y\in B} |f(y)-f_{B}|m(y)
\end{equation}
where the supremum is taken over all balls $B$ with radius $r(B)>0$ and $f_{B}$ denotes the mean value of $f$ on $B$ defined by \eqref{mean}. 
\begin{rem} \label{Nf=0}
For further use, observe that, if $f$ is a non-constant function on $\Gamma$, then $Nf(x)\neq 0$ for all $x\in \Gamma$. Indeed, if $Nf(x)=0$ for some $x\in \Gamma$, then $f(y)=f_B$ for all balls $B$ containing $x$. Thus, $f$ is constant on any ball containing $x$, therefore constant on $\Gamma$.
\end{rem}
The second maximal function we use is inspired by \cite{art} and \cite{maxhs}. Its definition involves estimates on the (discrete) divergence of test functions. More precisely, for all function $f$ on $\Gamma$, define, for all $x\in \Gamma$,
\begin{equation}\label{max gradiant f}
\mathcal{M}^{+}(f)(x)=\sup_{F}\left\vert \sum_{y\in \overset\circ{B}}f(y)(\delta F)(y)m(y)\right\vert,
\end{equation}
where the supremum is taken over all balls $B\subset\Gamma$ containing $x$ and all antisymmetric functions $F:E\rightarrow \R$ supported in $E_B$ and satisfying
\begin{equation}\label{les conditions sur F}
 |\vert F |\vert_{L^{\infty}(E_B)} \leq \frac{1}{V(B)},\ \ \ \ \ |\vert \delta F |\vert_{L^{\infty}(\overset\circ{B})} \leq \frac{1}{r(B)V(B)}.
\end{equation}
Define now, for $1\leq p\leq +\infty$,
$$
S^{1,p}(\Gamma):=\left\{f\in L^p(\Gamma); Nf\in L^p(\Gamma)\right\},
$$
equipped with the norm
$$
\left\Vert f\right\Vert_{S^{1,p}(\Gamma)}:=\left\Vert f\right\Vert_{L^p(\Gamma)}+\left\Vert Nf\right\Vert_{L^p(\Gamma)}.
$$
Consider also the $\dot{S}^{1,p}(\Gamma)$ space, made of functions $f$ on $\Gamma$ such that $Nf\in L^p(\Gamma)$, equipped with the semi-norm
$$
\left\Vert f\right\Vert_{\dot{S}^{1,p}(\Gamma)}:=\left\Vert Nf\right\Vert_{L^p(\Gamma)}.
$$
Define also
$$
E^{1,p}(\Gamma):=\left\{f\in L^p(\Gamma); \mathcal{M}^+f\in L^p(\Gamma)\right\},
$$
equipped with the norm
$$
\left\Vert f\right\Vert_{E^{1,p}(\Gamma)}:=\left\Vert f\right\Vert_{L^p(\Gamma)}+\left\Vert \mathcal{M}^+f\right\Vert_{L^p(\Gamma)},
$$
as well as its homogenous version. \par
\noindent Our first result is that, under \eqref{doubling}, \eqref{deltaalpha} and \eqref{p1}, the spaces $W^{1,p}(\Gamma)$, $S^{1,p}(\Gamma)$, $E^{1,p}(\Gamma)$ and $M^{1,p}(\Gamma)$, as well as their homogenous versions, coincide:
\begin{theo} \label{characsobolev}
Let $1<p\leq +\infty$. Assume that $\Gamma$ satisfies  \eqref{doubling}, \eqref{deltaalpha} and \eqref{p1}. Then:
\begin{itemize}
\item[$1.$]
$W^{1,p}(\Gamma)=S^{1,p}(\Gamma)=E^{1,p}(\Gamma)=M^{1,p}(\Gamma)$,
\item[$2.$]
$\dot{W}^{1,p}(\Gamma)=\dot{S}^{1,p}(\Gamma)=\dot{E}^{1,p}(\Gamma)=\dot{M}^{1,p}(\Gamma)$.
\end{itemize}
\end{theo}
\subsection{Characterization of Hardy-Sobolev spaces}
When $p=1$, as in the Euclidean case recalled in the introduction, the conclusion of Theorem \ref{characsobolev} does not hold. The following example is inspired by \cite{hajlasznew}, Example $3$. Take $\Gamma=\Z$ with its usual metric. Define, for all $x\in \Z$,
$$
f(x):=
\left\{
\begin{array}{ll}
\displaystyle \frac{x}{\left\vert x\right\vert \ln \left\vert x\right\vert} &\mbox{ if }\left\vert x\right\vert\geq 2,\\
0 & \mbox{ if }\left\vert x\right\vert \leq 1.
\end{array}
\right.
$$
Then $f\in \dot{W}^{1,1}(\Z)$. Indeed, for all $x\geq 2$, the mean-value theorem yields
$$
\left\vert f(x+1)-f(x)\right\vert=\left\vert \frac 1{\ln x}-\frac 1{\ln (x+1)}\right\vert\leq \frac 1{x\left(\ln x\right)^2}.
$$
As a consequence, for all $x\geq 3$,
\begin{equation} \label{gradf}
\left\vert \nabla f(x)\right\vert\leq \frac C{\left\vert x\right\vert \left(\ln \left\vert x\right\vert \right)^2}.
\end{equation}
Since $f$ is odd, \eqref{gradf} also holds for all $x\leq -3$. As a consequence, 
$$
\sum_{x\in \Z} \left\vert \nabla f(x)\right\vert<+\infty.
$$
Assume now that there exists a non-negative function $g \in L^1(\Z)$ such that 
$\vert f(x)-f(y)\vert \leq d(x,y) \left(g(x)+g(y)\right)$ for all $x,y\in \Z$. Then, for all $x\geq 3$,
$$
\left\vert f(x)-f(-x)\right\vert\leq 2x \left(g(x)+g(-x)\right).
$$
Since $f$ is odd, this means that, for all $x\geq 3$,
$$
\frac 1x\left\vert f(x)\right\vert\leq \left(g(x)+g(-x)\right).
$$
Therefore,
$$
2\sum_{\left\vert x\right\vert\geq 3} g(x)\geq \sum_{x\geq 3} \frac 1{x\ln x}=+\infty,
$$
which contradicts the fact that $g\in L^1(\Z)$.\par
\noindent The goal of this section is to give an endpoint version of Theorem \ref{characsobolev} when $p=1$. We will focus on the case of homogenous spaces. As it will turn out, asssuming \eqref{doubling} and $(P_1)$, one still has $\dot{M}^{1,1}(\Gamma)=\dot{S}^{1,1}(\Gamma)$.  Two extra characterizations of $\dot{M}^{1,1}(\Gamma)$ will be given: the first one is formulated in terms of ${\mathcal M}^+f$, the second one is an atomic decomposition. We first introduce these new descriptions.
\subsubsection{Maximal Hardy-Sobolev space}
It turns out that, as in the Euclidean case and in the context of Riemannian manifolds (see the introduction), Hardy-Sobolev spaces on $\Gamma$ can be defined by means of the functional ${\mathcal M}^{+}$. Let us first give a definition:
\begin{defi}\label{maximal Hardy-Sobolev space}
\textbf{(Maximal Hardy-Sobolev space)} 
\begin{itemize}
\item[$1.$]
We define the Hardy-Sobolev space $HS^{1}_{\max}(\Gamma)$ as follows:
\begin{equation}
HS^{1}_{\max}(\Gamma)=\{f\in L^{1}(\Gamma): \mathcal{M}^{+}f\in L^{1}(\Gamma)\}.
\end{equation}
This space is equipped with the norm 
\begin{equation}\label{norm hardy-sobolev max}
 |\vert f|\vert_{HS^{1}_{\max}(\Gamma)}:=\left\Vert f\right\Vert_{L^1(\Gamma)}+ \left\Vert \mathcal{M}^{+}f\right\Vert_{L^1(\Gamma)}.
\end{equation}
\item[$2.$]
The homogenous Hardy-Sobolev space $\dot{H}S^1_{\max}(\Gamma)$ is the space of all functions $f$ on $\Gamma$ such that ${\mathcal M}^{+}f\in L^1(\Gamma)$. It is equipped with the semi-norm
$$
 |\vert f|\vert_{\dot{H}S^{1}_{\max}(\Gamma)}:=\left\Vert \mathcal{M}^{+}f\right\Vert_{L^1(\Gamma)}.
$$
\end{itemize}
\end{defi}
\subsubsection{Atomic Hardy-Sobolev spaces}
\begin{defi}\label{Atom}
 For $1<t \leq +\infty$, define $t'$ by  $\frac{1}{t}+\frac{1}{t'}=1$. Say that a function $a$ on $\Gamma$ is a homogeneous
 Hardy-Sobolev $(1,t)-atom $ if
\begin{itemize}
\item[$1.$]
$a$ is supported in a ball $B$,
\item[$2.$]
$|\vert \nabla a |\vert_{t}\leq V(B)^{-\frac{1}{t'}}$,
\item[$3.$]
$\sum_{x\in \Gamma} a(x)m(x)=0$.
\end{itemize}
\end{defi}
If $f$ is a function on $\Gamma$, say that $f\in \dot{H}S^{1}_{t,ato}(\Gamma)$ if there exist a sequence $(\lambda_i)_{i\geq 1}\in l^1$ and a sequence of homogeneous Hardy-Sobolev $(1,t)$-atoms such that
\begin{equation}\label{decatomic}
 f=\sum\limits_{i}\lambda_{i} a_{i}.
\end{equation}
This space is equipped with the semi-norm
$$|\vert f |\vert_{\dot{H}S^{1}_{t,ato}(\Gamma)} = \inf \sum\limits_{i} |\lambda_{i}|$$
where the infimum is taken over all possible decompositions. \par
\noindent Notice that the convergence in \eqref{decatomic} is required to hold in $\dot{W}^{1,1}(\Gamma)$, which means that
$$
\lim_{k\rightarrow +\infty} \left\Vert \nabla \left(f-\sum_{j=0}^k \lambda_ja_j\right)\right\Vert_{L^1(\Gamma)}=0.
$$
The link between convergence in \eqref{decatomic} and pointwise convergence will be made explicit in Proposition \ref{pointwiseconv} below.\par
\noindent In the sequel, we will establish:
\begin{theo} \label{equal1}
Assume that \eqref{doubling}, \eqref{deltaalpha} and $(P_1)$ hold. Then $\dot{S}^{1,1}(\Gamma)=\dot{M}^{1,1}(\Gamma)=\dot{H}S^1_{\max}(\Gamma)=\dot{H}S^1_{t,ato}(\Gamma)$ for all $t\in (1,+\infty]$. In particular, $\dot{H}S^1_{t,ato}(\Gamma)$ does not depend on $t$.
\end{theo}
\begin{rem} \label{relax}
Assume that, in Definition \ref{Atom}, we replace condition $3$ by
\begin{itemize}
\item[$3^{\prime}$]
$\left\Vert a\right\Vert_{L^t(B)}\leq rV(B)^{-\frac{1}{t'}}$,
\end{itemize}
where $r$ is the radius of $B$, and we define $\dot{H}S^{1}_{t,ato}(\Gamma)$ as before, using this new type of atoms. Then, as the proof of Theorem \ref{equal1} will show (see Remark \ref{relax2} below), we obtain exactly the same $\dot{H}S^{1}_{t,ato}(\Gamma)$ space. This remark (inspired by ideas in \cite{badr2010atomic}) will turn out to be important for the study of Riesz transforms.
\end{rem}
\subsection{Interpolation}
As a consequence of the characterization of Hardy-Sobolev and Sobolev spaces through maximal functions, we establish an interpolation result between Hardy-Sobolev and Sobolev spaces:
\begin{theo} \label{interpol}
Let $1<q\leq +\infty$ and $\theta\in (0,1)$. Define $p$ such that $\frac 1p=(1-\theta)+\frac{\theta}q$. Then, for the complex interpolation method,
$$
\left[\dot{S}^{1,1}(\Gamma),\dot{W}^{1,q}(\Gamma)\right]_{\theta}=\dot{W}^{1,p}(\Gamma).
$$
\end{theo}
\subsection{Riesz transforms}
The Riesz transform in our context is the operator $R:=d(I-P)^{-1/2}$, which maps functions on $\Gamma$ to functions on $E$. The equality \eqref{bypartlaplacian} shows that $R$ is $L^2(Ð\Gamma)-L^2(E)$ bounded. For $1<p<+\infty$, the $L^p$-boundedness of $R$ was investigated in \cite{BR09} under various assumptions\footnote{Observe that the $L^p$-boundedness results of \cite{BR09} are stated for the operator $\nabla(I-P)^{-1/2}$, but \eqref{la relation entre diff f et grad f} shows at once that analogous conclusions hold for $d(I-P)^{-1/2}$.}. In particular, under \eqref{doubling} and the Poincar\'e inequality $(P_2)$, $R$ is $L^p(\Gamma)-L^p(E)$ bounded for all $1<p\leq 2$ (and even under weaker assumptions, see \cite{scand}). \par
\noindent For $p=1$, the Riesz transform is not $L^1(\Gamma)-L^1(E)$ bounded, but an endpoint version of the $L^p$-boundedness of $R$ for $1<p\leq 2$ was proved in \cite{pota}. This endpoint version involves the $H^1(\Gamma)$ atomic Hardy space on $\Gamma$, the definition of which we recall now. An atom in $H^1(\Gamma)$ is a function $a\in L^2(\Gamma)$, supported in a ball $B\subset \Gamma$ and satisfying
$$
\sum_{x\in \Gamma} a(x)m(x)=0\mbox{ and }\left\Vert a\right\Vert_{L^2(\Gamma)}\leq V(B)^{-1/2}.
$$
A function $f$ on $\Gamma$ is said to belong to $H^1(\Gamma)$ if and only if there exist a sequence $(\lambda_j)_{j\geq 1}\in l^1$ and a sequence of atoms $(a_j)_{j\geq 1}$ such that
$$
f=\sum_{j} \lambda_ja_j,
$$
where the series converges in $L^1(\Gamma)$. In this case, define
$$
\left\Vert f\right\Vert_{H^1(\Gamma)}:=\inf\sum_{j} \left\vert \lambda_j\right\vert,
$$
where, as usual, the infimum is taken over all possible decompositions of $f$. \par
\noindent Under \eqref{doubling} and $(P_2)$, the Riesz transform is $H^1(\Gamma)-L^1(E)$ bounded (\cite{pota}). This means that $(I-P)^{-1/2}$ is bounded from $H^1(\Gamma)$ to $\dot{W}^{1,1}(\Gamma)$. Here, under an extra assumption on the volume growth of balls of $\Gamma$, we prove that $(I-P)^{-1/2}$ maps continuously $H^1(\Gamma)$ into $\dot{S}^{1,1}(\Gamma)$:
\begin{theo} \label{Riesz}
Assume that $\Gamma$ satisfies \eqref{doubling} and $(P_2)$. Assume furthermore that there exist $C>0$ and $d\geq 1$ such that, for all $x\in \Gamma$ and all $1\leq r\leq s$,
\begin{equation} \label{reversevol}
\frac{V(x,r)}{V(x,s)}\leq C\left(\frac rs\right)^d.
\end{equation}
Then $(I-P)^{-1/2}$ is bounded from $H^1(\Gamma)$ into $\dot{S}^{1,1}(\Gamma)$.
\end{theo}
\begin{rem}
Under \eqref{doubling}, there exists $C^{\prime}>0$ such that, for all $x\in \Gamma$ and all $r\geq 1$,
$$
V(x,C^{\prime}r)\geq 2V(x,r)
$$
(see \cite{cougri}, Lemma 2.2). This implies that \eqref{reversevol} always holds with some $d>0$. In Theorem \ref{Riesz}, we assume furthermore that $d\geq 1$. This technical assumption seems to be required by our argument (see the proof of Theorem \ref{Riesz} in Section \ref{Rieszproofs} below), and could probably be removed. Note that assumption \eqref{reversevol} is satisfied when, for instance, $V(x,r)\sim r^d$ for some $d\geq 1$, which holds when $\Gamma$ is the Cayley graph of a group with polynomial volume growth.
\end{rem}
\section{Proofs of the characterizations of Sobolev spaces}
This section is devoted to the proof of Theorem \ref{characsobolev}. It will be convenient to use the following observation:
\begin{lem} \label{oscillNf}
For all functions $f$ on $\Gamma$, all $x\in \Gamma$ and all $r\geq 0$,
\begin{equation} \label{eqoscill}
\left\vert f(x)-f_{B(x,r)}\right\vert \leq CrNf(x).
\end{equation}
\end{lem}
{\bf Proof of Lemma \ref{oscillNf}: }  first, the conclusion is trivial when $0\leq r<1$, since in this case, $B(x,r)=\left\{x\right\}$ so that the left-hand side of \eqref{eqoscill} vanishes. Assume now that $r\geq 1$ and let $j\in \N$ be the integer such that $2^j\leq r<2^{j+1}$. Define $B:=B(x,2^{j+1})$ and, for all $-1\leq i\leq j+1$, $B_i=B(x,2^{i})$, so that $B=B_{j+1}$. Since $f(x)=f_{B\left(x,\frac 12\right)}$,
\begin{equation}\label{estimation fx-fB}
\begin{array}{lll}
\displaystyle \vert f(x)- f_{B}\vert&\leq&\displaystyle  \sum\limits_{i=-1}^{j}\left\vert f_{B_i}-f_{B_{i+1}}\right\vert\\
&\leq&\displaystyle \sum\limits_{i=-1}^{j} \frac{1}{V(B_{i})}
 \sum\limits_{y\in B_i} \vert f(y)-f_{B_{i+1}}\vert m(y)\\
 &\leq&\displaystyle C \sum\limits_{i=-1}^{j} \frac{r(B_{i+1})}{r(B_{i+1})V(B_{i+1})} \sum\limits_{y
 \in B_{i+1}} \vert f(y)-f_{B_{i+1}}\vert m(y)\\
 &\leq&\displaystyle C 2^jNf(x),
\end{array}
\end{equation}
where the third line uses \eqref{doubling}. Moreover, since $B(x,r)\subset B$,
\begin{equation} \label{term2}
\begin{array}{lll}
\displaystyle \left\vert f_{B(x,r)}-f_{B}\right\vert & \leq & \displaystyle \frac 1{V(x,r)}\sum_{y\in B(x,r)} \left\vert f(y)-f_B\right\vert m(y)\\
& \leq & \displaystyle C\frac 1{V(B)}\sum_{y\in B} \left\vert f(y)-f_B\right\vert m(y)\\
& \leq & \displaystyle C2^jNf(x),
\end{array}
\end{equation}
and the conjunction of \eqref{estimation fx-fB} and \eqref{term2} yields the conclusion (note that we used \eqref{doubling} again in the second line). \hfill\fin\par
\noindent As a corollary, one has (see also Lemma 3.6 in \cite{quasicont}):
\begin{prop} \label{fNf}
For all functions $f$ on $\Gamma$ and all $x,y\in \Gamma$,
$$
\left\vert f(x)-f(y)\right\vert\lesssim d(x,y)\left(Nf(x)+Nf(y)\right).
$$
\end{prop}
{\bf Proof: } let $x,y\in \Gamma$ with $x\neq y$ and $r:=d(x,y)$. Lemma \ref{oscillNf} yields
\begin{equation} \label{one}
\left\vert f(x)-f_{B(x,r)}\right\vert\leq CrNf(x).
\end{equation}
On the other hand, since $B(x,r)\subset B(y,2r)$, using Lemma \ref{oscillNf} again, one obtains
\begin{equation}\label{estimation fy-fB}
\begin{array}{lll}
\displaystyle \vert f(y)- f_{B(x,r)} \vert& \leq &\displaystyle \vert f(y) - f_{B(y,2r)} \vert + 
\vert f_{B(y,2r)}-f_{B(x,r)}\vert\\
&\leq& \displaystyle CrNf(y)+\frac{1}{V(x,r)}\sum\limits_{z\in B(x,r)} \left\vert f(z)-f_{B(y,2r)}
\right\vert m(z)\\
&\leq & \displaystyle C rNf(y) +C \frac{1}{V(y,2r)}\sum\limits_{z\in B(y,2r)} \left\vert f(z)-f_{B(y,2r)}\right\vert m(z)\\
&\leq&\displaystyle C r Nf(y).
\end{array}
\end{equation}
Thus, \eqref{one} and \eqref{estimation fy-fB} yield the desired result. \hfill\fin\par

\medskip

To establish that Sobolev spaces can also be characterized in terms of ${\mathcal M}^+f$, we have to solve the equation $\delta F=g$ in $L^{\infty}$ spaces
(see also \cite{maxhs}, Proposition 5.1 and \cite{dmrt} for the original ideas):
\begin{prop}\label{solution of delta F}
Assume that $\Gamma$ satisfies $(D)$ and $(P_{1})$. Let $B$ a ball of $\Gamma$ with $r(B)\geq 1$ and $g\in L_{0}^{\infty}(B)$. Then, there exists $F\in L^{\infty}(E_{B})$ such that
$\delta F= g\mbox{ in }\overset{\circ}{B}$ and
\begin{equation}\label{estimation of sobolev}
 |\vert F |\vert_{L^{\infty}(E_{B})} \lesssim r(B) |\vert g|\vert_{L^{\infty}(B)}.
\end{equation}
\end{prop}
\textbf{Proof: } let $B$ be a ball and $g\in L_{0}^{\infty}(B)$. Consider 
$$\mathcal{S}=\{V\in L^{1}(E_{B}): \exists f\in L^{1}(\Gamma)\mbox{ supported in }\overset{\circ}{B},\ V=df\mbox{ in }E_B\}.$$
We consider $\mathcal{S}$ as subspace of $L^{1}(E_{B})$ equipped with the norm 
$$|\vert V |\vert_{L^{1}(E_{B})} = \sum\limits_{\gamma\in E_{B}}|V(\gamma)|\mu_{\gamma}$$
(see Section \ref{measures}). Define a linear functional on $\mathcal{S}$ by 
$$L(V):=\sum\limits_{x\in B} g(x)f(x) m(x)\mbox{ if } V=df\in \mathcal{S}.$$
Observe that $L$ is well defined since $\sum\limits_{x\in B}g(x)m(x)=0$ and it is plain to see that, if $df_1=df_2$ in $E_B$, then $f_1-f_2$ is constant on $B$. From $(P_{1})$ and using the support condition on $f$, we derive 
$$
\begin{array}{lll}
 |L(V)|&\leq & \displaystyle \sum\limits_{x\in B}\left\vert g(x)\right\vert |f(x)-f_{B}|m(x)\\
&\leq& \displaystyle Cr(B) |\vert g |\vert_{L^{\infty}(B)}\sum\limits_{x\in B } \nabla f(x)m(x)\\
& \leq & \displaystyle Cr(B)\left\Vert g\right\Vert_{L^{\infty}(B)} \sum_{x\in B} \left(\sum_{y\sim x} p(x,y)\left\vert f(y)-f(x)\right\vert\right)m(x)\\
& = & \displaystyle Cr(B)\left\Vert g\right\Vert_{L^{\infty}(B)} \sum_{x\in B} \sum_{y\sim x}\left\vert f(y)-f(x)\right\vert \mu_{xy}\\
& = & \displaystyle Cr(B)\left\Vert g\right\Vert_{L^{\infty}(B)} \sum_{x\sim y,\ x\in B,\ y\in B} \left\vert f(y)-f(x)\right\vert \mu_{xy}\\
&=& \displaystyle Cr(B)|\vert g |\vert_{L^{\infty}(B)} |\vert V |\vert_{L^{1}(E_{B})}.
\end{array}
$$
The Hahn-Banach theorem shows that $L$ can be extended to a bounded linear functional on $L^{1}(E_B)$ with norm not greater than 
$Cr(B)|\vert g |\vert_{\infty}$. Thus, there exists $F\in L^{\infty}(E_B)$ such that, for all $V\in L^1(E_B)$,
$$
L(V)=\sum_{\gamma\in E_B} F(\gamma)V(\gamma)\mu_{\gamma}.
$$
In particular, for all $f\in L^1(B)$ vanishing on $\partial B$, \eqref{IPP} yields\footnote{Observe that $F$ and $df$ are square integrable on $E_B$ since $E_B$ is a finite set.}
$$
\sum_{x\in B} g(x)f(x)m(x)=L(df)=\sum_{\gamma\in E_B} F(\gamma)df(\gamma)\mu_{\gamma}=-\sum_{x\in \overset{\circ}{B}} \delta F(x)f(x)m(x),
$$
which ensures that $-\delta F=g$ in $\overset{\circ}{B}$ with
$$
\left\Vert F\right\Vert_{L^{\infty}}\leq Cr(B)|\vert g |\vert_{\infty}.
$$
\hfill\fin\par
\noindent A consequence of Proposition \ref{solution of delta F} , which will also be useful in the proof of Theorem \ref{characsobolev}, is:
\begin{prop} \label{MN}
For all functions $f$ on $\Gamma$:
\begin{itemize}
\item[$1.$]
$$
{\mathcal M}^+f\sim Nf,
$$
\item[$2.$]
$$
\nabla f \lesssim Nf.
$$
\end{itemize}
\end{prop}
{\bf Proof of ${\mathcal M}^+f\lesssim Nf$: } let $x\in \Gamma$. Take $F$ 
as in the definition of ${\mathcal M}^+f$, associated to a ball $B$ containing $x$. Then \eqref{IPP}, applied with the function $f$ equal to $1$ in $\overset{\circ}B$ and to $0$ on $\partial B$, shows that
$\sum\limits_{y\in \overset{\circ}{B}} \left(\delta F \right)(y)m(y)=0$ so we can write 
$$\left|\sum\limits_{y\in \overset{\circ}{B}}f(y)(\delta F)(y)m(y)\right|=\left|\sum\limits_{y\in \overset{\circ}{B}}\left(f(y)-f_{B}\right)(\delta F)(y)m(y)\right|.$$
Thus, \eqref{les conditions sur F} yields 
$$\left|\sum\limits_{y\in \overset{\circ}{B}}f(y)(\delta F)(y)m(y)\right| \leq \frac{1}{r(B)V(B)} \sum\limits_{y\in \overset{\circ}{B}} \left|f(y)-f_{B}\right|m(y)$$
$$\lesssim Nf(x).$$
Taking the supremum over all such $F$, we get 
$$
{\mathcal M}^+f(x)\lesssim Nf(x).
$$
\noindent {\bf Proof of $Nf\lesssim {\mathcal M}^+f$: } let $x\in \Gamma$ and $B=B(x_B,r(B))$ a ball containing $x$. We may and do assume that $r_B\geq 1$, otherwise
$$
\sum_{y\in B} \left\vert f(y)-f_B\right\vert m(y)=0.
$$
Define $\widetilde{B}:=B(x_B,r(B)+1)$, so that $B\subset\overset{\circ}{\widetilde{B}}$. If $g\in L^{\infty}_{0}(B)$ with  $|\vert g |\vert_{\infty}\leq 1$, extend $g$ by $0$ outside $B$ and solve $\delta F=g$ in $\overset{\circ}{\widetilde{B}}$ with $F\in L^{\infty}\left(E_{\widetilde{B}}\right)$ satisfying \eqref{estimation of sobolev}. Extend $F$ by $0$ outside $E_{\widetilde{B}}$ . Then, setting
$$\widetilde{F}:=\frac{F}{Cr\left(\widetilde{B}\right)V(\widetilde{B})},$$ one has
$$
\begin{array}{lll}
\displaystyle \frac 1{r(B)V(B)} \left\vert \sum_{y\in B} f(y)g(y)m(y)\right\vert & = & \displaystyle  \frac 1{r(B)V(B)} \left|\sum\limits_{y\in \overset{\circ}{\widetilde{B}}} f(y) g(y)m(y)\right|\\
&= &\displaystyle \frac 1{r(B)V(B)} \left|\sum\limits_{y\in \overset{\circ}{\widetilde{B}}} f(y)(\delta F)(y)m(y)\right|\\
&=& \displaystyle C\frac{r(\widetilde{B})V(\widetilde{B})}{r(B)V(B)} \left|\sum\limits_{y\in \overset{\circ}{\widetilde{B}}}f(y) (\delta \widetilde{F})(y)m(y)\right|\\
& \leq & \displaystyle C{\mathcal M}^+f(x),
\end{array}
$$
where the last line follows from \eqref{doubling} and the fact that $\widetilde{F}$ satisfies \eqref{les conditions sur F} . 
Taking the supremum on the left hand side over all balls containing $x$, we get $Nf(x)\leq C\mathcal{M}^{+}f(x)$. This inequality concludes the proof of $1.$\par
\noindent {\bf Proof of $\nabla f\lesssim Nf$: } let $x\in \Gamma$. Fix $y\sim x$, set $B:=B(x,2)$ and define the function $F$ on $E$ in the following way: $F(x,y)=\frac 1{m(x)}$, $F(y,x)=-\frac 1{m(x)}$ and $F(u,v)=0$ whenever $(u,v)\neq (x,y)$ and $(u,v)\neq (y,x)$. Notice that $\delta F$ is supported in $\overset{\circ}{B}$ and
$$
\left\Vert F\right\Vert_{L^{\infty}(E_B)}\lesssim \frac 1{V(B)}\mbox{ and }\left\Vert \delta F\right\Vert_{L^{\infty}(\overset{\circ}{B})}\lesssim \frac 1{r(B)V(B)}.
$$
This and item $1$ of Proposition \ref{MN} yield 
$$
\left\vert \langle df,F\rangle\right\vert=\left\vert \langle f,\delta F\rangle\right\vert \lesssim {\mathcal M}^+f(x)\lesssim Nf(x).
$$
But
$$
\langle df,F\rangle=2(f(y)-f(x))\frac{\mu_{xy}}{m(x)}=2p(x,y)(f(y)-f(x)),
$$
which shows that
$$
p(x,y)\left\vert f(y)-f(x)\right\vert\lesssim Nf(x)
$$
for all $y\sim x$. The definition of $\nabla f$ then yields the desired result. \hfill\fin\par

\medskip

\noindent {\bf Proof of Theorem \ref{characsobolev}: } we write it for homogenous spaces, the inhomogeneous case being an immediate consequence. First, assertion $1$ in Proposition \ref{MN} gives at once that $\dot{E}^{1,p}(\Gamma)=\dot{S}^{1,p}(\Gamma)$.\par
\noindent Assume now that $f\in \dot{W}^{1,p}(\Gamma)$ and let $x\in \Gamma$. Since $(P_p)$ holds, there exists $q<p$ such that $(P_q)$ is still valid (see Remark \ref{improvepoinc}). For all balls $B\ni x$, $(P_q)$ yields
$$
\frac 1{V(B)}\sum_{y\in B} \left\vert f(y)-f_B\right\vert m(y)\leq Cr(B) \left(\frac 1{V(B)}\sum_{y\in B} \left\vert \nabla f(y)\right\vert^qm(y)\right)^{\frac 1q},
$$
so that, taking the supremum over $B$,
$$
Nf(x)\leq C\left({\mathcal M}_{HL}\left\vert\nabla f\right\vert^q\right)^{\frac 1q}(x),
$$
where ${\mathcal M}_{HL}$ stands for the Hardy-Littlewood maximal function, given by
$$
{\mathcal M}_{HL}f(x):=\sup_{B\ni x} \frac 1{V(B)}\sum_{y\in B} \left\vert f(y)\right\vert m(y),
$$
where, again, the supremum is taken over all balls $B$ containing $x$. Since $\nabla f\in L^p(\Gamma)$ and ${\mathcal M}_{HL}$ is $L^{\frac pq}(\Gamma)$-bounded (this is because \eqref{doubling} holds and $\frac pq>1$), one has
$$
\left(\sum_{x\in \Gamma} \left\vert Nf(x)\right\vert^pm(x)\right)^{\frac 1p} \leq C\left\Vert \left({\mathcal M}_{HL}\left\vert\nabla f\right\vert^q\right)\right\Vert_{L^{\frac pq}(\Gamma)}^{\frac 1q} \leq C\left\Vert \nabla f\right\Vert_{L^p(\Gamma)},
$$
which shows that $Nf\in L^p(\Gamma)$. One therefore has $f\in \dot{S}^{1,p}(\Gamma)$ and $\left\Vert f\right\Vert_{\dot{S}^{1,p}(\Gamma)}\leq C\left\Vert f\right\Vert_{\dot{W}^{1,p}(\Gamma)}$. \par
\noindent  Take now $f\in \dot{S}^{1,p}(\Gamma)$. Since $Nf\in L^p(\Gamma)$, Proposition \ref{fNf} shows that $f\in \dot{M}^{1,p}(\Gamma)$ and $\left\Vert f\right\Vert_{\dot{M}^{1,p}(\Gamma)}\leq C\left\Vert f\right\Vert_{\dot{S}^{1,p}(\Gamma)}$. \par
\noindent Assume finally that $f\in \dot{M}^{1,p}(\Gamma)$ and let $g\in L^p(\Gamma)$ given by \eqref{Hajlasz space} and satisfying $\left\Vert g\right\Vert_{L^p(\Gamma)}\leq 2\left\Vert f\right\Vert_{\dot{M}^{1,p}(\Gamma)}$. Define, for all $x\in \Gamma$, $h(x):=\sum_{y\sim x} \left(g(y)+g(x)\right)$. Then $h\in L^p(\Gamma)$ and $\left\Vert h\right\Vert_{L^p(\Gamma)}\leq C\left\Vert g\right\Vert_{L^p(\Gamma)}$. Indeed, observing that, whenever $x\sim y$, $m(x)\leq Cm(y)$ (this is an immediate consequence of \eqref{doubling}), and using the fact that any point in $\Gamma$ has at most $N$ neighbours, one obtains
$$
\begin{array}{lll}
\displaystyle \sum_{x\in \Gamma} h(x)^pm(x) & \leq & \displaystyle C\sum_{x\sim y} \left(g(x)^p+g(y)^p\right)m(x)\\
& \leq & \displaystyle C\sum_{x\in \Gamma} g(x)^pm(x)+C\sum_{y\in \Gamma} g(y)^pm(y)\\
& = & \displaystyle C\left\Vert g\right\Vert_{L^p(\Gamma)}^p.
\end{array}
$$
Now, let $x\in \Gamma$. By \eqref{Hajlasz space} and the fact that $0\leq p(x,y)\leq 1$ for all $x,y\in \Gamma$,
$$
\nabla f(x)\leq C\sum_{y\sim x} \left\vert f(y)-f(x)\right\vert\leq C\sum_{y\sim x} \left(g(x)+g(y)\right)=Ch(x),
$$
so that $\nabla f\in L^p(\Gamma)$ and $\left\Vert f\right\Vert_{\dot{W}^{1,p}(\Gamma)}\leq C\left\Vert f\right\Vert_{\dot{M}^{1,p}(\Gamma)}$. This completes the proof. \hfill\fin
\section{The Calder\'on-Zygmund decomposition for Hardy-Sobolev spaces}
The present section is devoted to the proof of the Calder\'on-Zygmund decomposition for Hardy-Sobolev spaces on graphs. The corresponding decomposition on Riemannian manifolds was established in \cite{badr2010atomic}. Recall that analogous Calder\'on-Zygmund decompositions for classical Sobolev spaces were proved in \cite{auscher2005riesz} on Riemannian manifolds and \cite{BR09} on graphs.
\begin{prop}\label{Calderon-Zygmund decomposition}
[Calder\'{o}n-Zygmund decomposition for Hardy-Sobolev spaces] Let $\Gamma$ satisfy $(D)$ and $(P_1)$. 
Let $f\in \dot{S}^{1,1}(\Gamma), \frac{s}{s+1}< q <1$ and $\alpha>0$. Then one can find a collection of balls $\{B_{i}\}_{i\in I}$, functions $b_{i}\in W^{1,1}(\Gamma)$ 
and a function $g\in \dot{W}^{1,\infty}(\Gamma)$ such that the 
following properties hold:
$$f=g+\sum\limits_{i}b_{i},$$
\begin{equation}\label{gradiant g} 
 |\nabla g(x)|\leq C \alpha \mbox{ for all } x\in \Gamma,
\end{equation}
\begin{equation}\label{support bi dans Bi}
 supp \ \ b_{i} \subset B_{i},  \Vert  b_i \Vert_1 \leq C \alpha r_iV(B_i), \ \ \Vert \nabla b_i
 \Vert_q \leq C \alpha V(B_i)^{1/q}
\end{equation}
\begin{equation}\label{sum les volums de balls}
 \sum\limits_{i}V(B_{i})\leq \frac{C}{\alpha} \sum\limits_{x\in B_i}(Nf)(x) m(x)
\end{equation}
and 
\begin{equation}\label{sum sur Bi}
 \sum\limits_{i}\chi_{B_{i}}\leq K,
\end{equation}
where, for all $i$, $r_i$ is the radius of $B_i$, and $C$ and $K$ only depend on $q,p$ and on the constants in $(D)$ and $(P_{1})$.
\end{prop}
{\bf Proof: }Êthe proof of Proposition \ref{Calderon-Zygmund decomposition} follows the main lines of the one of Proposition 3.3 in \cite{badr2010atomic}, with adaptations due to the discrete context.\par
\noindent Let $f\in \dot{S}^{1,1}(\Gamma)$ and $\alpha> 0$. Define 
$$\Omega:=\left\{x\in \Gamma;\ \mathcal{M}_{HL,q}(Nf)(x)> \frac {\alpha}C \right\},$$
where $C$ is the implicit constant in item $2$ of Proposition \ref{MN} and ${\mathcal M}_{HL,q}$ is defined by
\begin{equation} \label{defmq}
\mathcal{M}_{HL,q}(g)(y):=\left({\mathcal M}_{HL} \left\vert g\right\vert^q)\right)^{1/q}.
\end{equation}
Let $F:=\Gamma\setminus \Omega$.\par
\noindent A consequence of item $2$ in Proposition \ref{MN} is that
\begin{equation} \label{nablafF}
\nabla f(x)\leq CNf(x)\leq C{\mathcal M}_{HL,q}(Nf)(x)\leq \alpha\mbox{ for all }x\in F.
\end{equation}
If $\Omega=\emptyset $, then set 
$$f=g \mbox{ and }b_{i}=0\mbox{ for all }i,$$
so that \eqref{gradiant g} is satisfied by \eqref{nablafF}, and all the other required properties are clearly satisfied. \par
\noindent From now on, assume that $\Omega\neq \emptyset$. First,
\begin{equation} \label{Nf is bounded}
        \begin{array}{lll}
m(\Omega)&\leq& \displaystyle \frac{C}{\alpha} \sum\limits_{x\in \Gamma} \mathcal{M}_{HL,q}(Nf)(x)m(x)\\
         &=& \displaystyle \frac{C}{\alpha} \sum\limits_{x\in \Gamma} \left(\mathcal{M}_{HL}(|Nf|^{q})(x)\right)^{\frac{1}{q}}m(x) \\
         & \leq & \displaystyle \frac{C}{\alpha}\sum\limits_{x\in \Gamma} Nf(x)m(x)< 
         \infty,
         \end{array}
\end{equation}
where, in the last line, we used the fact the $\mathcal{M}_{HL}$ is $L^{1/q}(\Gamma)$-bounded since
 $q<1$ and $Nf\in L^1(\Gamma)$. In particular $\Omega 
\neq \Gamma$ as $m(\Gamma )= + \infty$ (see Remark \ref{mart}). \par
\noindent {\bf Definition of the balls $B_i$: } since $\Omega$ is a strict subset of $\Gamma$, let $\{\underline{B_{i}}\}_{i}$
 be a Whitney decomposition of $\Omega$ (see \cite{coifmanweiss}). More precisely, the $\underline{B_{i}}$ are pairwise disjoint, and there exist two constants $C_{2}> C_{1}>1,$ only depending on the metric,
such that 
\begin{itemize}
 \item $\Omega=\cup_{i} B_{i}$ with $B_{i}= C_{1}\underline{B_{i}},$ and the balls $B_{i}$ have the bounded overlap property,
\item $r_{i}=r(B_{i})= \frac{1}{2} d(x_{i},F)$ where $x_{i}$ is the center of $B_{i}$,
\item each ball $\overline{B_{i}}=C_{2}B_{i}$ intersects $F$ (one can take $C_{2}=4C_{1}$). 
\end{itemize}
For $x \in \Omega$, define $I_{x}:=\{i: x\in B_{i}\}$. As already seen in \cite{BR09}, there exists $K$ such that $\sharp I_{x}\leq K,$ and moreover, for all $i,k\in I_{x}$, $\frac 13 r_{i}\leq r_{k}\leq 3r_{i}$ and $B_{i}\subset 7B_{k}$. The bounded overlap property yields \eqref{sum sur Bi} and implies
\begin{equation} \label{sumvolbi}
\sum_{i} V(B_i)\lesssim m(\Omega).
\end{equation}
Then, \eqref{sum les volums de balls} follows from \eqref{sum sur Bi} and \eqref{Nf is bounded}. \par
\noindent The following observation will be used several times: for all $i$,
\begin{equation} \label{majorNf}
\left(\frac 1{V(C_2B_i)} \sum\limits_{x\in C_2B_{i}}|Nf(x)|^{q}m(x)\right)^{\frac{1}{q}}\leq C\alpha V(B_i).
\end{equation}
Indeed, the left-hand side of \eqref{majorNf} is bounded by ${\mathcal M}_{HL,q}(Nf)(y)$ for some $y\in C_2B_i\cap F$, which yields the result. \par
\noindent {\bf Definition of the functions $b_{i}$: } following the construction in Section 5 of \cite{BR09}, pick up a partition of unity $\{\chi_{i}\}_{i}$ of $\Omega$ 
subordinated to the covering $\{B_{i}\}_{i}$. Each $\chi_{i}$ is a Lipschitz function supported in $B_{i}$ with $0\leq \chi_{i}\leq 1,\ \vert|\nabla \chi_{i}|\vert_{\infty}\leq \frac{C}{r_{i}}$ and 
$\sum\limits_{i\in I} \chi_i (x)={\bf 1}_{\Omega}$  for all $x\in\Gamma$. Moreover, $\nabla \chi_i$ is supported in $C_3 B_i\subset \Omega$ with $C_3< 2.$  We set $b_{i}:=(f-f_{B_i}) \chi_{i}$, so that $supp\ b_i \subset B_i$.\par
\noindent {\bf Estimate of $\left\Vert b_i\right\Vert_{L^1(\Gamma)}$: } the Sobolev-Poincar\'{e} inequality \eqref{sobpoincm11} applied with $g=Nf$ (recall that $Nf\in L^1(\Gamma)$ and the pair ($f,Nf)$ satisfies \eqref{Hajlasz space} by Proposition \ref{fNf}) and $\lambda=C_2$, as well as \eqref{majorNf}, yield
\begin{equation}\label{better estimate bi}
\begin{array}{lll}
\displaystyle  \Vert b_{i}\Vert_{1}&\leq & \displaystyle  \sum\limits_{x\in B_{i}} |f(x)-f_{B_{i}}|
m(x)\\
&\leq& \displaystyle C r_{i} \left(\frac 1{V(C_2B_i)} \sum\limits_{x\in C_2B_{i}}|Nf(x)|^{q}m(x)\right)^{\frac{1}{q}}
V(B_{i})\\
& \leq & \displaystyle C r_{i} \alpha V(B_{i}).
\end{array}
\end{equation}
{\bf Proof of $\nabla b_i\in L^{1}(\Gamma)$: } since 
$$\nabla b_i (x)=\nabla\left((f-f_{B_i})\chi_i\right)(x)\leq \left(\max\limits_{y\sim x}\chi_i(y)\right)\nabla 
f(x)+\left\vert f(x)-f_{B_i}\right\vert \nabla \chi_i(x)$$
and $\chi_i\leq 1$ on $\Gamma$, using \eqref{majorNf} again, one obtains  
\begin{equation}\label{estimation bi}
\begin{array}{lll}
 \displaystyle |\vert \nabla b_{i}|\vert_{1}&\leq& \displaystyle \sum\limits_{x\in C_3B_{i}} \vert 
 f(x)-f_{B_i}\vert \vert \nabla \chi_{i}(x)\vert m(x)+ \sum\limits_{x\in C_3B_i}\vert \nabla f(x)\vert m
 (x)\\
&\leq& \displaystyle C \alpha V(B_{i})+\sum\limits_{x\in C_3B_i} |\nabla f(x)|m(x)<+\infty.
\end{array}
\end{equation}
{\bf Estimate of $\left\Vert \nabla b_i\right\Vert_{L^q(\Gamma)}$: } using item $2$ in Proposition \ref{MN}, \eqref{sobpoincm11} with $g=Nf$ (and H\"older) and \eqref{majorNf}, we obtain:
\begin{equation}\label{estimate bi in the sobolev space}
\begin{array}{lll}
\Vert \nabla b_{i}\Vert_{q}^q& \leq & \displaystyle C \left( \sum\limits_{x\in C_3B_i} \vert \nabla 
f(x)\vert^q m(x) + \sum\limits_{x\in C_3B_i} \left\vert f(x)-f_{B_i}\right\vert^q\vert\nabla\chi_i(x)\vert^q 
m(x)\right)\\
&\leq & \displaystyle C \sum\limits_{x\in C_2B_i} \vert Nf(x)\vert^qm(x)+C \frac{C^q}{r_i^q} 
r_i^q \left(\sum\limits_{x\in C_2B_i} \vert Nf(x)\vert^q m(x)\right)\\
&\leq &C \alpha^q V(B_{i}).
\end{array}
\end{equation}
Thus \eqref{support bi dans Bi} is proved.\par

\medskip

\noindent {\bf Definition of $g$: }set now $g=f-\sum\limits_{i} b_{i}.$ Since the sum is locally finite on $\Omega$, $g$ is well-defined on $\Gamma$ and $g=f$ on $F$.\par
\noindent {\bf Estimate of $\left\vert \nabla g\right\vert$: } since $\sum\limits_{i\in I} \chi_i(x)=1$ for all $x
\in \Omega$, one has 
\begin{equation*}
g=f\chi_F + \sum\limits_{i\in I} f_{B_i} \chi_i
\end{equation*}
where $\chi_F$ denotes the characteristic function of $F$. We will need the following 
lemma:
\begin{lem}\label{gugv}
There exists $C>0$ such that, for all $j\in I$, all $u\in F\cap 4B_j$ and all $v\in B_j,$
$$\vert g(u)-g(v)\vert \leq C \alpha d(u,v).$$
\end{lem}
Let us admit the conclusion of Lemma \ref{gugv} and complete the proof of \eqref{gradiant g}. It is enough to check that $\vert g(x) - g(y)\vert \leq C \alpha$ 
for all $x\sim y\in \Gamma$. Three situations may occur:
\begin{enumerate}
\item Assume first that $x,y\in\Omega$. Let $j\in I$ such that $x\in B_j$. Since $\chi_F(x)= \chi_F(y)=0$ and $\sum_i \chi_i=1$ on $\Gamma$, it follows that
$$g(y)-g(x)=\sum\limits_{i\in I} \left(f_{B_i}-f_{B_j}\right)\left(\chi_i(y)-\chi_i(x)\right),$$
so that $\vert g(y)- g(x)\vert \leq C\sum\limits_{i\in I}\vert f_{B_i}-f_{B_j} \vert 
\nabla\chi_i(x):= h(x)$. \par
\noindent We claim that $\vert h(x)\vert \leq C\alpha$, which will end the proof in this case. 
Let $i\in I$ be such that $\nabla \chi_i(x)\neq 0$, so that $d(x,B_i)\leq 1$, hence $r_i \leq 3 r_j+1\leq 4r_j$ and
$B_i\subset 10B_j$. An application of \eqref{sobpoincm11} with $g=Nf$ and of \eqref{majorNf} yields
\begin{equation}\label{estimation fBi moins fB10j}
\begin{array}{lll}
\displaystyle \vert  f_{B_i} - f_{10B_j} \vert & \leq&  \displaystyle \frac{1}{V(B_i)}\sum
\limits_{y\in B_i} \vert f(y)- f_{10B_j} \vert m(y)\\
&\leq& \displaystyle \frac{C}{V(B_j)}\sum\limits_{y\in 10B_j} \vert f(y)- f_{10 B_j}\vert m
(y)\\
&\leq &  \displaystyle C r_j \left( \frac{1}{V(10B_j)} \sum\limits_{y\in 10B_j}\vert Nf(y) \vert^q m(y)\right)^{1/q}\\
&\leq&  \displaystyle C r_j \alpha.
\end{array}
\end{equation}
Analogously  $\vert f_{10B_j}-f_{B_j}\vert \leq C r_j \alpha. $ Hence
\begin{equation}\label{estimation hx}
\begin{array}{lll}
\displaystyle \vert h(x)\vert & =& \displaystyle \left\vert \sum\limits_{i\in I; x\in 2B_i}(f_
{B_i}-f_{B_j})\nabla \chi_i(x)\right\vert\\
&\leq& \displaystyle C  \sum\limits_{i\in I; x\in 2B_i} \vert f_{B_i } -f_{B_j}\vert r_{i}^{-1}\\
&\leq& C K \alpha.
\end{array}
\end{equation}
\item Assume now that $x\in \overset\circ{F}$, so that $y\in F$. In this case $\vert g(y)-g(x)\vert = \vert f(x)- f(y) 
\vert \leq C\nabla f(x)\leq C \alpha$ by \eqref{nablafF}.
\item Assume finally that $x\in \partial F$.
\subitem i.\ If $y\in F$, as already seen, $\vert g(y)-g(x)\vert= \vert f(x)-f(y)\vert \leq C \nabla f(x)
 \leq C\alpha$ by \eqref{nablafF}.
 \subitem ii.\ Assume finally that $y\in \Omega$. There exists  $j\in I$ such that  
 $y\in B_j$. Since $x\sim y,$ one has $x\in 4B_j, $ Lemma \ref{gugv} therefore yields
 \begin{equation*}
 \vert g(x)- g(y)\vert \leq C \alpha d(x,y) \leq C \alpha.
 \end{equation*}
\end{enumerate}
The case when $x\in \Omega$ and $y\in F$ is contained in Case 3.ii by symmetry, since $y\in \partial F$.  Thus the proof of Proposition \ref{Calderon-Zygmund decomposition} is complete. \hfill\fin\par
\medskip
\noindent {\bf Proof of Lemma  \ref{gugv}: } it is analogous to the one of Lemma 5.1 in \cite{BR09}.  The only difference is that one uses \eqref{sobpoincm11} instead of the Poincar\'e inequality applied in \cite{BR09}. \hfill\fin \par

\section{Proofs of the characterization of Hardy-Sobolev spaces} \label{pfs}
We now turn to the proof of Theorem \ref{equal1}. Let us explain the strategy. We first establish that $\dot{S}^{1,1}(\Gamma)=\dot{M}^{1,1}(\Gamma)$. The inclusion $\dot{S}^{1,1}(\Gamma)\subset \dot{M}^{1,1}(\Gamma)$ is proved exactly in the same way as the corresponding inclusion in Theorem \ref{characsobolev}. The converse is more involved, since the Hardy-Littlewood maximal function is not $L^1(\Gamma)$-bounded, and the proof relies on a Sobolev-Poincar\'e inequality. \par
\noindent The identity $\dot{S}^{1,1}(\Gamma)=\dot{H}S^1_{\max}(\Gamma)$ is an immediate consequence of item $1$ in Proposition \ref{MN}.\par
\noindent Finally, we check that $\dot{S}^{1,1}(\Gamma)=\dot{H}S^1_{ato}(\Gamma)$, using the Sobolev-Poincar\'e inequality again, as well as an adapted Calder\'on-Zygmund decomposition.
\subsection{Sharp maximal characterization of $\dot{M}^{1,1}(\Gamma)$}
A straightforward consequence of Proposition \ref{fNf} is that $\dot{S}^{1,1}(\Gamma)\subset \dot{M}^{1,1}(\Gamma)$. \par
\noindent The proof of the converse inclusion relies, as the proof of Theorem 3 in \cite{kin2007}, on a Sobolev-Poincar\'e inequality (\cite{hajlasz}, theorem 8.7) :
\begin{theo} \label{sobpoinc}
Let $p\in \left[\frac{s}{s+1},s\right)$\footnote{where $s$ is given by \eqref{defs}.}, $B\subset \Gamma$ be a ball with radius $r$, $f\in \dot{M}^{1,p}(B)$ and $g\in L^{p}(B)$ such that $(f,g)$ satisfies \eqref{Hajlasz space} in $B$ (see Remark \ref{m1pleq1}). Then $(f,g)$ satisfies the following Sobolev-Poincar\'{e} inquality: for all $\lambda>1$, there is a constant $C> 0$ only depending on the constant in \eqref{doubling} and $\lambda$ such that
\begin{equation}\label{la relation entre u et g}
 \left(\frac{1}{V(B)}\sum\limits_{x\in B}|f(x)-f_{B}|^{p^*}m(x)\right)^{\frac{1}{p^*}} \leq C r \left(\frac{1}{V(\lambda B)}
\sum\limits_{x\in \lambda B} g(x)^{p} m(x)\right)^{\frac{1}{p}}
\end{equation}
where $p^*:=\frac{sp}{s-p}.$
\end{theo}
An easy consequence of Theorem \ref{sobpoinc} is that, for all functions $f\in \dot{M}^{1,1}(\Gamma)$, all $q\in \left[\frac{s}{s+1},s\right)$, all balls $B\subset \Gamma$ of radius $r$ and all $\lambda>1$,
\begin{equation} \label{sobpoincm11}
\frac 1{V(B)} \sum_{x\in B} \left\vert f(x)-f_B\right\vert m(x)\leq Cr\left(\frac{1}{V(\lambda B)}
\sum\limits_{x\in \lambda B} g(x)^{q} m(x)\right)^{\frac{1}{q}}
\end{equation}
whenever $(f,g)$ satisfies \eqref{Hajlasz space}. Indeed, it is enough to observe that $g\in L^{\frac{s}{s+1}}(\lambda B)$, apply Theorem \ref{sobpoinc} with $p=\frac s{s+1}$, since $p^{\ast}=1$ and use H\"older inequality. \par
\medskip
\noindent Take now $f\in \dot{M}^{1,1}(\Gamma)$, $q\in \left[\frac{s}{s+1},1\right)$ and $g$ such that \eqref{Hajlasz space} and \eqref{sobpoincm11} hold and $\left\Vert g\right\Vert_{L^1(\Gamma)}\leq 2\left\Vert f\right\Vert_{\dot{M}^{1,1}(\Gamma)}$. The inequality \eqref{sobpoincm11} yields
$$
Nf(y)\lesssim {\mathcal M}_{HL,q}g(y)
$$
for all $y\in \Gamma$, where ${\mathcal M}_{HL,q}$ was defined by \eqref{defmq}. Since $1/q>1$, the Hardy-Littlewood maximal function is $L^{1/q}(\Gamma)$-bounded, which implies that
$$
\left\Vert Nf\right\Vert_{L^1(\Gamma)}\lesssim \left\Vert g\right\Vert_{L^1(\Gamma)}\lesssim \left\Vert f\right\Vert_{\dot{M}^{1,1}(\Gamma)}.
$$
This ends the proof of the inclusion $\dot{M}^{1,1}(\Gamma)\subset \dot{S}^{1,1}(\Gamma)$. \hfill\fin
\subsection{Maximal characterization}
The identity $\dot{S}^{1,1}(\Gamma)=\dot{H}S^1_{\max}(\Gamma)$ is an immediate consequence of item $1$ in Proposition \ref{MN}.
\subsection{Atomic decomposition}
We prove now that $\dot{H}S^{1}_{t,ato}(\Gamma)=\dot{S}^{1,1}(\Gamma)$ for all $t\in (1,+\infty]$.
\subsubsection{$\dot{H}S^{1}_{t,ato}(\Gamma)\subset \dot{S}^{1}_{1}(\Gamma)$}
For the proof of this inclusion, we have to clarify the link between convergence in $\dot{H}S^1_{t,ato}(\Gamma)$ and pointwise convergence:
\begin{prop} \label{pointwiseconv}
Let $f\in \dot{H}S^1_{t,ato}(\Gamma)$ and write
$$
f=\sum_{j} \lambda_ja_j,
$$
where $\sum_{j} \left\vert \lambda_j\right\vert<+\infty$, for all $j$, $a_j$ is a homogeneous Hardy-Sobolev $(1,t)$-atom and the series converges in $\dot{W}^{1,1}(\Gamma)$. Then, for all $k$, there exists $c_k\in \R$ such that, for all $x\in \Gamma$,
$$
f(x)=\lim_{k\rightarrow +\infty} \sum_{j=0}^k \lambda_ja_j(x)-c_k.
$$
\end{prop}
The proof follows from:
\begin{lem} \label{pointl1}
Let $(h_k)_{k\geq 1}\in \dot{W}^{1,1}(\Gamma)$. If $\lim_{k\rightarrow +\infty} \left\Vert \nabla h_k\right\Vert_{L^1(\Gamma)}=0$, then, for all $k\geq 1$, there exists $c_k\in \R$ such that
$$
\lim_{k\rightarrow +\infty} h_k(x)-c_k=0.
$$
\end{lem}
{\bf Proof of Lemma \ref{pointl1}: }Êassume first that there exists $x_0\in \Gamma$ such that $h_k(x_0)=0$ for all $k\geq 1$. Then, for all $x\in \Gamma$, $\lim_{k\rightarrow +\infty} h_k(x)=0$. Indeed, the very definition of $\nabla h_k$ implies that, for all $x,y\in \Gamma$ with $x\sim y$, $\lim_{k\rightarrow +\infty} (h_k(x)-h_k(y))=0$. The conclusion then readily follows for all $j\geq 1$ and for all $x\in B(x_0,j)$ by induction on $j$. \par
\noindent In the general case, fix $x_0\in \Gamma$ and define $g_k(x):=h_k(x)-h_k(x_0)$ for all $k\geq 1$ and all $x\in \Gamma$. What we have just seen means that $\lim_{k\rightarrow +\infty} g_k(x)=0$, which yields the desired conclusion with $c_k:=h_k(x_0)$. \hfill\fin\par
\noindent {\bf Proof of Proposition \ref{pointwiseconv}: } it is an immediate consequence of Lemma \ref{pointl1} applied with $h_k:=f-\sum_{j=0}^k \lambda_ja_j$. \hfill\fin
\begin{prop}\label{HSA in M}
Assume that $\Gamma$ satisfies $(D)$ and $(P_{1})$. Let $t\in (1,+\infty]$.
\begin{itemize}
\item[$1.$]
Let $a$ be a homogeneous $(1,t)$ atom. Then $a\in \dot{S}_{1}^{1}(\Gamma)$ with $|\vert a |\vert_{S^{1}_{1}}\leq C$.
\item[$2.$]
One has $\dot{H}S^{1}_{t,ato}(\Gamma)\subset \dot{S}^{1}_{1}(\Gamma)$ and there exists $C>0$ such that, for all $f\in  \dot{H}S^{1}_{t,ato}(\Gamma)$,
$$|\vert f |\vert_{\dot{S}^{1}_{1}(\Gamma)}\leq C |\vert f |\vert_{\dot{H}S^{1}_{ato}(\Gamma)}.$$
\end{itemize}
\end{prop}
\textbf{Proof: } for $1$, let $a$ be a homogeneous $(1,t)$ atom supported in ball $B=B(x,r).$ We want to prove that $Na \in L^{1}(\Gamma)$ and that $\left\Vert Na\right\Vert_{L^1(\Gamma)}\leq C$. For all $y\in \Gamma$, and all balls $B^{\prime}\ni y$, $(P_{1})$ yields:
$$
\begin{array}{lll}
 \displaystyle \frac{1}{r(B^{\prime})V(B^{\prime})} \sum\limits_{z\in B^{\prime}}|a(z)-a_{B^{\prime}}|m(z)
      &\leq& \displaystyle  \frac{C}{V(B^{\prime})}\sum\limits_{z\in B^{\prime}} \nabla a(z) m(z)\\
      & \leq  & \displaystyle {\mathcal M}_{HL}(\nabla a)(y),
\end{array}
$$
so that
\begin{equation} \label{estimNa}
Na(y)\lesssim {\mathcal M}_{HL}(\nabla a)(y).
\end{equation}
As a consequence, 
\begin{equation} \label{Na2B}
\begin{array}{lll}
\displaystyle \sum_{y\in B(x,4r)} Na(y)m(y) & \leq & \displaystyle CV(x,4r)^{1/t^{\prime}} \left(\sum_{y\in B(x,4r)} \left({\mathcal M}_{HL}(\nabla a)(y)\right)^tm(y)\right)^{1/t}\\
& \leq & \displaystyle CV(x,4r)^{1/t^{\prime}} \left\Vert \nabla a\right\Vert_{L^t(\Gamma)}\\
& \leq & C,
\end{array}
\end{equation}
where the first line follows from H\"older and \eqref{estimNa}, the second one from the $L^t$-boundedness of the Hardy-Littlewood maximal function and the last one from the doubling property and the second item in Definition \ref{Atom}.\par
\noindent Let $k\geq 2$ and $y\in B(x,2^{k+1}r)\setminus B(x,2^kr)$. Consider an arbitrary ball $B^{\prime}$ containing $y$. One has
$$
\begin{array}{lll}
\displaystyle \frac 1{r(B^{\prime})V(B^{\prime})} \sum_{z\in B^{\prime}} \left\vert a(z)-a_{B^{\prime}}\right\vert m(y) &= & \displaystyle  \frac 1{r(B^{\prime})V(B^{\prime})} \sum_{z\in B^{\prime}\cap B} \left\vert a(z)-a_{B^{\prime}}\right\vert m(z) \\
& +  & \displaystyle \frac 1{r(B^{\prime})V(B^{\prime})} \sum_{z\in B^{\prime}\setminus B} \left\vert a_{B^{\prime}}\right\vert m(z) \\
& \leq & \displaystyle \frac 3{r(B^{\prime})V(B^{\prime})} \sum_{z\in B^{\prime}\cap B} \left\vert a(z)\right\vert m(z).
\end{array}
$$
It is easily checked that, if $B^{\prime}\cap B\neq\emptyset$, then $r(B^{\prime})>2^{k-1}r$ and \eqref{doubling} yields $V(x,2^{k+1}r)\leq CV(B^{\prime})$. As a consequence of this observation and $(P_1)$ (remember that $a_B=0$),
$$
\begin{array}{lll}
Na(y) & \leq & \displaystyle \frac C{2^{k-1}rV(2^{k+1}B)} \sum_{z\in B} \left\vert a(z)\right\vert m(z)\\
& \leq  & \displaystyle \frac C{2^{k-1}V(2^{k+1}B)} \sum_{z\in B} \left\vert \nabla a(z)\right\vert m(z)\\
& \leq & \displaystyle \frac C{2^{k-1}V(2^{k+1}B)}.
\end{array}
$$
It follows that
\begin{equation} \label{Naoutside2B}
\begin{array}{lll}
\displaystyle \sum_{y\notin B(x,4r)} Na(y)m(y)&= & \displaystyle \sum_{k\geq 2} \sum_{y\in B(x,2^{k+1}r)\setminus B(x,2^kr)} Na(y)m(y)\\
& \leq & \displaystyle \sum_{k\geq 2} \frac C{2^{k-1}V(2^{k+1}B)} V(2^{k+1}B)\\
& \leq &  C.
\end{array}
\end{equation}
Gathering \eqref{Na2B} and \eqref{Naoutside2B}, one obtains $\left\Vert Na\right\Vert_{L^1(\Gamma)}\leq C$.\par
\noindent Now, for assertion $2$ in Proposition \ref{HSA in M}, if $f\in \dot{H}S^{1}_{t,ato}(\Gamma)$, take an atomic decomposition of $f:$ 
$f=\sum\limits_{i} \lambda_{i}a_{i}$ where each $a_{i}$ is an atom and $\sum\limits_{i} |\lambda_{i}|\leq 2\left\Vert f\right\Vert_{\dot{H}S^{1}_{t,ato}(\Gamma)}$. By Proposition \ref{pointwiseconv}, pick up a sequence $(c_k)_{k\geq 1}\in \R$ such that, for all $x\in \Gamma$,
$$
f(x)=\lim_{k\rightarrow +\infty} \sum_{j=0}^k \lambda_ja_j(x)-c_k=\lim_{k\rightarrow +\infty} f_k(x)-c_k
$$
where, for all $k$, $f_k:=\sum_{j=0}^k \lambda_ja_j$. \par
\noindent Let $x\in \Gamma$ and $B$ be a ball containing $x$. Observe that
$$
f_B=\frac 1{V(B)}\sum_{y\in B} f(y)m(y)=\lim_{k\rightarrow +\infty} \frac 1{V(B)} \sum_{y\in B} \left(f_k(y)-c_k\right)m(y)=\lim_{k\rightarrow +\infty} \left((f_k)_B-c_k\right).
$$
As a consequence,
$$
\frac 1{V(B)}\sum_{y\in B} \left\vert f(y)-f_B\right\vert m(y)=\lim_{k\rightarrow +\infty} \frac 1{V(B)} \sum_{y\in B} \left\vert f_k(y)-(f_k)_B\right\vert m(y).
$$
For all $k\geq 1$,
$$
\sum_{y\in B} \left\vert f_k(y)-(f_k)_B\right\vert m(y) \leq  \sum_{j=0}^k \left\vert \lambda_j\right\vert\sum_{y\in B}  \left\vert a_j(y)-(a_j)_B\right\vert m(y),
$$
so that
$$
\frac 1{r(B)V(B)}\sum_{y\in B} \left\vert f(y)-f_B\right\vert m(y)\leq \sum_{j=0}^{+\infty} \left\vert \lambda_j\right\vert Na_j(x).
$$
Since $\left\Vert Na_j\right\Vert_{L^1(\Gamma)}\leq C$ and $\sum_j \left\vert \lambda_j\right\vert\leq 
2\left\Vert f\right\Vert_{\dot{H}S^{1}_{t,ato}(\Gamma)}$, Proposition \ref{HSA in M} is proved.
\hfill\fin\par
\begin{rem} \label{relax2}
Observe that, in the above argument, if condition $3$ in Definition \ref{Atom} is replaced by condition $3^{\prime}$ in Remark \ref{relax}, then the previous computation is still valid, since one has, using H\"older,
$$
\begin{array}{lll}
Na(y) & \leq & \displaystyle \frac C{2^{k-1}rV(2^{k+1}B)} \sum_{z\in B} \left\vert a(z)\right\vert m(z)\\
& \leq  & \displaystyle \frac C{2^{k-1}rV(2^{k+1}B)} \left\Vert a\right\Vert_{L^t(B)} V(B)^{1/t^{\prime}}\\
& \leq & \displaystyle \frac C{2^{k-1}V(2^{k+1}B)}.
\end{array}
$$
\end{rem}

\subsubsection{$\dot{S}^{1}_{1}(\Gamma)\subset \dot{H}S^{1}_{q^{\ast},ato}(\Gamma)$}

\bigskip

\noindent The proof of the inclusion $\dot{S}^{1,1}(\Gamma)\subset \dot{H}S^{1}_{ato}(\Gamma)$ relies on the 
Calder\'{o}n-Zygmund decomposition for functions in $\dot{S}^{1,1}(\Gamma)$ given by Proposition \ref{Calderon-Zygmund decomposition}:

\begin{prop}\label{decomposition f}
Let $\Gamma$ satisfying \eqref{doubling} and $(P_1)$. Let $f\in \dot{S}^{1,1}(\Gamma)$. Then for all $\frac{s}{s+1}< q < 1, q^{*}=\frac{sq}{s-q},$ there is a sequence of
 $(1,q^{*})$ Hardy-Sobolev atoms $\{a_{j}\}_{j},$ and a sequence of scalars $\{\lambda _{j}\}_{j}\in l^1$ such that 
$$f=\sum\limits_{j} \lambda_{j}a_{j} \ \ \ \mbox{in} \ \dot{W}^{1,1}(\Gamma), \ \ \mbox{and} \ \ \ \sum|\lambda_{j}|\leq 
C_{q} |\vert f|\vert_{\dot{S}^{1,1}(\Gamma)}.$$
Consequently, $\dot{S}^{1,1}(\Gamma)\subset \dot{H}S^{1}_{q^{*},ato}(\Gamma)$ with $|\vert f |\vert_{\dot{H}S^{1}_{q^{*},ato}(\Gamma)}\leq C_{q} |\vert f |\vert_{\dot{S}^{1,1}(\Gamma)}.$
\end{prop}
{\bf Proof: } the proof is analogous to the one of Proposition 3.4 in \cite{badr2010atomic}, which deals with the case of Riemannian manifolds, and is also inspired by the proof of the atomic decomposition for Hardy spaces in \cite{E.M.Stein93}, section III.2.3. We may and do assume that $f$ is not constant on $\Gamma$, otherwise one can take $a_j=0$ for all $j$. \par
\noindent Let $f\in \dot{S}^{1,1}(\Gamma)$. For every $j\in \Z^{\ast}$, we take the Calder\'{o}n-Zygmund decomposition for $f$ at level $\alpha=2^j$ given by Proposition \ref{Calderon-Zygmund decomposition}. Then $$f=g^j+\sum_i b^j_i$$ with $b^j_i, g^j$ satisfying the properties of Proposition \ref{Calderon-Zygmund decomposition}. We first claim
\begin{equation}\label{decomposition of f }
f=\sum^\infty_{-\infty} (g^{j+1}- g^j),
\end{equation}
where the series converges in $\dot{W}^{1,1}(\Gamma)$. \par
\noindent To see this, observe first that $g^j \rightarrow f$ in $\dot{W}^{1,1}(\Gamma)$ as $j \rightarrow +\infty$. Indeed, since the sum is locally finite we can write, using \eqref{estimation bi}, \eqref{sumvolbi} and the facts that $C_3B^j_i\subset \Omega$ and that the $C_3B^j_i$ have the bounded overlap property,
\begin{equation} \label{j+infty}
\begin{array}{lll}
\displaystyle \left\Vert \nabla (g^j - f)  \right\Vert_{L^1(\Gamma)} &=& \displaystyle \left\Vert \nabla 
\left(\sum_i b_i^j\right)\right\Vert_{L^1(\Gamma)} \leq \sum_i \left\Vert \nabla b_i^j \right\Vert_{L^1(\Gamma)} \\
&\leq &\displaystyle C2^j m(\Omega^j)+ C \sum_{x\in \Omega_j} \vert\nabla f(x) \vert m
(x)\\
&:=& \displaystyle I_j+II_j,
\end{array}
\end{equation}
where $\Omega^j:=\left\{x\in \Gamma , \mathcal{M}_{HL,q}(Nf)(x)> \frac{2^j}C\right\}$. Observe that $\Omega^{j+1} 
\subset \Omega^j$ for all $j\in \Z$.\par
\noindent Observe that
\begin{equation} \label{sumomegaj}
\sum_{j\in \Z} 2^jm(\Omega^j) \lesssim \int_0^{+\infty} m\left(\left\{x\in \Gamma;\ {\mathcal M}_{HL,q}(Nf)(x)>t\right\}\right)dt=\left\Vert {\mathcal M}_{HL,q}(Nf)\right\Vert_{L^1(\Gamma)}<+\infty.
\end{equation}
This implies that, when $j\rightarrow +\infty$, $I_j \rightarrow 0$. 
Since $\nabla f\in L^1(\Gamma)$ and $m(\Omega_j)\rightarrow 0$ when $j\rightarrow +\infty$, one has $II_j \rightarrow 0$ when $j\rightarrow +\infty$. Thus, \eqref{j+infty} shows that
$$
\lim_{j\rightarrow +\infty} g^j=f\mbox{ in }\dot{W}^{1,1}(\Gamma).
$$
\noindent Next, when $j\rightarrow - \infty$, we want to show $\left\Vert \nabla g_j 
\right\Vert_{L^1(\Gamma)} \rightarrow 0$. If $F^j:=\Gamma\setminus \Omega^j$, an immediate consequence of Remark \ref{Nf=0} is that, since $f$ is not constant on $\Gamma$,
\begin{equation} \label{empty}
\bigcap_{j\in \Z} F^j=\emptyset.
\end{equation}
Write
$$
\begin{array}{lll}
\displaystyle \left\Vert \nabla g^j 
\right\Vert_{L^1(\Gamma)} & \lesssim & \displaystyle \sum_{x\sim y,\ x,y\in F^j} \left\vert g^j(x)-g^j(y)\right\vert m(x)\\
& + & \displaystyle  \sum_{x\sim y,\ x,y\in \Omega^j} \left\vert g^j(x)-g^j(y)\right\vert m(x)\\
& + & \displaystyle  \sum_{x\sim y,\ x\in F^j,\ y\in \Omega^j} \left\vert g^j(x)-g^j(y)\right\vert m(x)\\
& := & A_j+B_j+C_j.
\end{array}
$$
If $x\sim y$ with $x\in F^j$ and $y\in F^j$, $\left\vert g^j(x)-gj(y)\right\vert=\left\vert f(x)-f(y)\right\vert$, so that
$$
A_j\lesssim \sum_{x\in F^j} \nabla f(x)m(x),
$$
which implies that $A_j\rightarrow 0$ when $j\rightarrow -\infty$, since $\nabla f\in L^1(\Gamma)$ and \eqref{empty} holds. \par
\noindent Moreover,
$$
B_j\lesssim \sum_{x\in \Omega_j} \nabla g^j(x) m(x)\lesssim 2^jm(\Omega^j),
$$
and this quantity goes to $0$ when $j\rightarrow -\infty$ by \eqref{sumomegaj}. \par
\noindent Finally, if $x\sim y$ with $x\in F^j$ and $y\in \Omega^j$, $\left\vert g^j(x)-g^j(y)\right\vert\lesssim \nabla g^j(y)$ and, since $m(x)$ and $m(y)$ are comparable when $x\sim y$, one has
$$
C_j\lesssim \sum_{y\in \Omega^j} \nabla g^j(y)m(y)
$$
which goes to $0$ when $j\rightarrow -\infty$. 
This ends the proof of \eqref{decomposition of f }.\par

\medskip

\noindent Introduce a partition of unity $(\chi^j_k)_{k}$ subordinated to balls $B^j_k$ corresponding to $\Omega^j$ as in the proof of Proposition \ref{Calderon-Zygmund decomposition}. We will need two observations:
\begin{lem} \label{observ}
\begin{itemize}
\item[$1.$]
For all $j,k,l$, if there exist $x\in B^j_k$ and $y\in B^{j+1}_l$ with $x\sim y$, then 
\begin{equation} \label{comparadius}
r^{j+1}_l\leq 4r_k^j.
\end{equation}
\item[$2.$]
There exists $C>0$ such that, for all $j$,
\begin{equation} \label{overlap}
\sum_k {\bf 1}_{2B^j_k}\leq C.
\end{equation}
\end{itemize}
\end{lem}
We postpone the proof of Lemma \ref{observ} and end up the proof of the atomic decomposition of $f$.\par
\noindent Set $g^{j+1} - g^j:=l^j$
and decompose $l_j$ as $\displaystyle l^j=\sum\limits_k l^j_k$ with 
\begin{equation}\label{ljk}
l^j_k:=(f-d^j_k)\chi^j_k - \sum_l (f-d^{j+1}_l) \chi_l^{j+1}\chi^j_k+\sum_l c_{k,l}^j \chi_l^{j+1},
\end{equation}
where, for all $j,k$,
$$
d^j_k:=\frac 1{\sum_{y} \chi^j_k(y)m(y)} \sum_{y} f(y)\chi^j_k(y)m(y),
$$
and
$$
\displaystyle c_{k,l}^j:= \frac{1}{\sum_{y\in B^{j+1}_l} \chi^{j+1}_l(y)m(y)}\sum_{x\in B_l^{j+1}} \left(f(x)-d^{j+1}_l\right)\chi_l^{j+1}(x) \chi_k^j(x)m(x).
$$
First, the identity $\displaystyle l^j=\sum\limits_k l^j_k$ holds by definition of $g^j$ and $g^{j+1}$ and since $\sum\limits_k \chi_k^j=1$ on the support of $\chi_l^{j+1}$ and, for all $l$, $\displaystyle \sum\limits_k c_{k,l}^j=0.$ \par
\noindent We now claim that, up to a constant, $2^{-j}V(B^j_k)^{-1}l^j_k$ is a homogeneous Hardy-Sobolev $(1,q^{\ast})$ atom.
Indeed, the cancellation condition
$$\sum_{x\in \Gamma} l^j_k(x) m(x)=0$$
for all $k$ follows from the fact that $\displaystyle \sum\limits_{x\in \Gamma}\left(f(x)-d^j_k\right)\chi_k^j(x) m(x)=0$ and 
 the definition of $c_{k,l}^j$, which immediately gives, for all $l$, $\displaystyle \sum\limits_{x\in \Gamma}  \left(\left(f(x)-d^{j+1}_l\right)\chi_l^{j+1}(x) \chi _k^j(x)- c_{k,l}^j \chi_l^{j+1}(x)\right)m(x)=0.$
 A consequence of \eqref{comparadius} is that $l^j_k$ is supported in the ball $9B^j_k$, therefore $\nabla l^j_k$ is supported in $18B^j_k$. \par
\noindent Let us now prove that 
\begin{equation}Ê\label{gradljk}
\left\Vert \nabla 
l^j_k \right\Vert_{L^{q^*}(\Gamma)}\lesssim 2^jV(B^j_k)^{1/q^{\ast}}.
\end{equation}
Let $x,y\in \Gamma$ such that $x\sim y$. Write
\begin{equation}
\begin{array}{lll}
\displaystyle l^j_k(y)-l^j_k(x) & = & \displaystyle \left((f(y)-f(x))\chi^j_k(y)-\sum_l (f(y)-f(x))\chi^{j+1}_l(y)\chi^j_k(y)\right) \\
& + & \displaystyle \left(f(x)-d^j_k\right) \left(\chi^j_k(y)-\chi^j_k(x)\right)\\
& - & \displaystyle \sum_l \left(f(x)-d^{j+1}_l\right)\left(\chi^{j+1}_l(y)\chi^j_k(y)-\chi^{j+1}_l(x)\chi^j_k(x)\right)\\
& + & \displaystyle \sum_l c_{k,l}\left(\chi^{j+1}_l(y)-\chi^{j+1}_l(x)\right)\\
& := & \Delta_1(x,y)+\Delta_2(x,y)+\Delta_3(x,y)+\Delta_4(x,y).
\end{array}
\end{equation}
Let us estimate $\Delta_i(x,y)$ for $1\leq i\leq 4$.\par
\noindent {\bf Estimate of $\Delta_1$: } compute
$$
\Delta_1(x,y)=(f(y)-f(x))\chi^j_k(y)\left(1-{\bf 1}_{\Omega^{j+1}}(y)\right).
$$
As a consequence, if $\Delta_1(x,y)\neq 0$, one has $y\in B^j_k\cap (\Omega^j\setminus \Omega^{j+1})$, so that $x\in 2B^j_k$. By item $2$ in Proposition \ref{MN}, one has $\nabla f(y)\leq C2^j$, so that $\left\vert f(y)-f(x)\right\vert \leq C2^j$. As a consequence, for all $x\in \Gamma$,
$$
\sum_{y\sim x} \left\vert \Delta_1(x,y)\right\vert^{q^{\ast}}\leq C2^{jq^{\ast}}.
$$
Therefore, by \eqref{doubling},
\begin{equation} \label{estimdelta1}
\sum_{x\in 2B^j_k} \sum_{y\sim x} \left\vert \Delta_1(x,y)\right\vert^{q^{\ast}}m(x)\leq C2^{jq^{\ast}} V(B^j_k).
\end{equation}
{\bf Estimate of $\Delta_2$: }Êobserve first that if $\Delta_2(x,y)\neq 0$, then $y\in B^j_k$ or $x\in B^j_k$, so that $x\in 2B^j_k$. Since $\nabla \chi^j_k\leq \frac C{r^j_k}$ on $\Gamma$, one has, for all $x\in \Gamma$,
$$
\sum_{y\sim x} \left\vert \Delta_2(x,y)\right\vert^{q^{\ast}}\leq \frac C{\left(r^j_k\right)^{q^{\ast}}} \left\vert f(x)-d^j_k\right\vert^{q^{\ast}} .
$$
As a consequence, 
$$
\sum_{x\in 2B^j_k} \sum_{y\sim x} \left\vert \Delta_2(x,y)\right\vert^{q^{\ast}} m(x) \leq   \frac C{\left(r^j_k\right)^{q^{\ast}}} \sum_{x\in 2B^j_k} \left\vert f(x)-d^j_k\right\vert^{q^{\ast}}m(x).
$$
But
$$
\left\Vert f-d^j_k\right\Vert_{L^{q^{\ast}}(2B^j_k)} \leq  \left\Vert f-f_{B^j_k}\right\Vert_{L^{q^{\ast}}(2B^j_k)}+\left\vert d^j_k-f_{B^j_k}\right\vert V^{1/q^{\ast}}(2B^j_k),
$$
and
$$
\begin{array}{llll}
\displaystyle V^{1/q^{\ast}}(2B^j_k) \left\vert d^j_k-f_{B^j_k}\right\vert & = & \displaystyle V^{1/q^{\ast}}(2B^j_k)\left\vert \frac 1{\sum_y \chi^j_k(y)m(y)} \sum_z \left(f(z)-f_{B^j_k}\right)\chi^j_k(z)m(z)\right\vert\\
& \leq & \displaystyle C\left(\frac {V(B^j_k)}{\sum_y \chi^j_k(y)m(y)} \sum_{z\in B^j_k} \left\vert f(z)-f_{B^j_k}\right\vert^{q^{\ast}}m(z)\right)^{1/q^{\ast}}\\
& \leq & \displaystyle C\left(\sum_{z\in B^j_k} \left\vert f(z)-f_{B^j_k}\right\vert^{q^{\ast}}m(z)\right)^{1/q^{\ast}}.
\end{array}
$$
Thus,
\begin{equation} \label{f-djk}
\sum_{x\in 2B^j_k} \left\vert f(x)-d^j_k\right\vert^{q^{\ast}}m(x)\leq C\sum_{z\in 2B^j_k} \left\vert f(z)-f_{B^j_k}\right\vert^{q^{\ast}}m(z).
\end{equation}
Therefore, by Theorem \ref{sobpoinc} and \eqref{majorNf},
$$
\begin{array}{lll}
\displaystyle \sum_{x\in 2B^j_k} \sum_{y\sim x} \left\vert \Delta_2(x,y)\right\vert^{q^{\ast}} m(x) & \leq  & \displaystyle \frac C{\left(r^j_k\right)^{q^{\ast}}} \sum_{x\in 2B^j_k} \left\vert f(x)-f_{B^j_k}\right\vert^{q^{\ast}}m(x) \\
& \leq & \displaystyle CV(B^j_k) \left(\frac 1{V(4C_2B^j_k)}\sum_{x\in 4C_2B^j_k} Nf(x)^qm(x)\right)^{\frac{q^{\ast}}q}\\
& \leq & \displaystyle CV(B^j_k)2^{jq^{\ast}}.
\end{array}
$$
{\bf Estimate of $\Delta_3(x,y)$: } first,
$$
\begin{array}{lll}
\displaystyle -\Delta_3(x,y) & = & \displaystyle \sum_l\left(f(x)-f_{B^{j+1}_l}\right)\chi^j_k(y)\left(\chi^{j+1}_l(y)-\chi^{j+1}_l(x)\right)\\
& + & \displaystyle \sum_l\left(f(x)-f_{B^{j+1}_l}\right)\chi^{j+1}_l(x)\left(\chi^j_k(y)-\chi^j_k(x)\right)\\
& = & \Delta_3^1(x,y)+\Delta_3^2(x,y).
\end{array}
$$
For $\Delta_3^1(x,y)$, notice that the sum may be computed over the $l\in I^j(x)$, where
$$
I^j(x):=\left\{l;\mbox{ there exists }y\sim x\mbox{ such that }y\in B^j_k\mbox{ and }x\mbox{ or }y\mbox{ belong to }B^{j+1}_l\right\}.
$$
For $l\in I^j(x)$, $x\in 2B^j_k\cap 2B^{j+1}_l$ and $r^{j+1}_l\leq 4r^j_k$ by Lemma \ref{observ}. Since $\left\vert \chi^{j+1}_l(y)-\chi^{j+1}_l(x)\right\vert\leq \frac C{r^{j+1}_l}$, one has, for all $x\in \Gamma$,
$$
\sum_{y\sim x} \left\vert \Delta_3^1(x,y)\right\vert^{q^{\ast}}\leq \sum_{l\in I^{j}(x)} \frac C{\left(r^{j+1}_l\right)^{q^{\ast}}} \left\vert f(x)-f_{B^{j+1}_l}\right\vert^{q^{\ast}}.
$$
Notice that, by item $2$ in Lemma \ref{observ}, $\sharp I^{j}(x)\leq C$. 
It follows that
$$
\begin{array}{lll}
\displaystyle \sum_{x\in 2B^j_k\cap 2B^{j+1}_l} \sum_{y\sim x} \left\vert \Delta_3^1(x,y)\right\vert^{q^{\ast}} m(x)& \leq  & \displaystyle \sum_{x\in 2B^j_k\cap 2B^{j+1}_l} \sum_{l\in I^{j}(x)} \frac C{\left(r^{j+1}_l\right)^{q^{\ast}}} \left\vert f(x)-f_{B^{j+1}_l}\right\vert^{q^{\ast}}m(x)\\
& = & \displaystyle C\sum_{l} \frac 1{\left(r^{j+1}_l\right)^{q^{\ast}}} \sum_{x\in 2B^j_k\cap 2B^{j+1}_l,\ l\in I^j(x)} \left\vert f(x)-f_{B^{j+1}_l}\right\vert^{q^{\ast}}m(x)\\
& \leq & \displaystyle  C\sum_{l;\ B^{j+1}_l\subset CB^j_k}V(CB^{j+1}_l) \left(\frac 1{V(4CB^{j+1}_l)} \sum_{x\in 4CB^{j+1}_l} Nf(x)^qm(x)\right)^{\frac{q^{\ast}}q}\\
& \leq & \displaystyle C\sum_{l;\ B^{j+1}_l\subset CB^j_k)}V(CB^{j+1}_l) 2^{(j+1)q^{\ast}}\\
& \leq & \displaystyle CV(CB^j_k)2^{jq^{\ast}}.
\end{array}.
$$
In this computation, we used the fact that, for $l\in I^{j}(x)$, one has $r^{j+1}_l\leq 4r^j_k$ and, since $2B^{j+1}_l\cap 2B^j_k\neq \emptyset$, $B^{j+1}_l\subset CB^j_k$.\par
\noindent For $\Delta_3^2(x,y)$, arguing similarly, the sum may be restricted to the $l\in J^j(x)$ where
$$
J^j(x):=\left\{l;\ x\in B^{j+1}_l\mbox{ and there exists }y\sim x\mbox{ such that }y\in B^j_k\mbox{ or }x\in B^j_k\right\}.
$$
For $l\in J^j(x)$, $x\in B^{j+1}_l\cap 2B^j_k$ and $r^{j+1}_l\leq 4r^j_k$. Again, $\sharp J^{j}(x)\leq C$. Arguing as before, one obtains
$$
\sum_{y\sim x} \left\vert \Delta_3^2(x,y)\right\vert^{q^{\ast}}\leq \sum_{l\in J^{j}(x)} \frac C{\left(r^{j}_k\right)^{q^{\ast}}} \left\vert f(x)-f_{B^{j+1}_l}\right\vert^{q^{\ast}}.
$$
As a consequence,
$$
\begin{array}{lll}
\displaystyle \sum_{x\in 2B^j_k\cap B^{j+1}_l} \sum_{y\sim x} \left\vert \Delta_3^2(x,y)\right\vert^{q^{\ast}} m(x)& \leq  & \displaystyle \sum_{x\in 2B^j_k\cap B^{j+1}_l} \sum_{l\in J^{j}(x)} \frac C{\left(r^{j}_k\right)^{q^{\ast}}} \left\vert f(x)-f_{B^{j+1}_l}\right\vert^{q^{\ast}}m(x)\\
& \leq & \displaystyle  \sum_{x\in 2B^j_k\cap B^{j+1}_l} \sum_{l\in J^{j}(x)} \frac C{\left(r^{j+1}_l\right)^{q^{\ast}}} \left\vert f(x)-f_{B^{j+1}_l}\right\vert^{q^{\ast}}m(x)\\
& \leq & \displaystyle C\sum_{l;\ B^{j+1}_l\subset CB^j_k }V(CB^{j+1}_l) \left(\frac 1{V(4CB^{j+1}_l)} \sum_{x\in 4CB^{j+1}_l} Nf(x)^qm(x)\right)^{\frac{q^{\ast}}q}\\
& \leq & \displaystyle C\sum_{l;\ B^{j+1}_l\subset CB^j_k} V(CB^{j+1}_l) 2^{(j+1)q^{\ast}}\\
& \leq & \displaystyle CV(CB^j_k)2^{jq^{\ast}}.
\end{array}.
$$
\noindent {\bf Estimate of $\Delta_4$: }Ê
note first that $c_{k,l}^j=0$ when 
 $B_k^j \cap B_l^{j+1}=\emptyset $ and $\vert c_{k,l}^j\vert\leq C 2^j r_l^{j+1}$ thanks 
 to \eqref{better estimate bi}. As a consequence, 
 $\left\vert c_{k,l}^j \left(\chi_l^{j+1}(y)-\chi^{j+1}_l(x)\right)\right\vert \leq C2^j$ for every $l$. It follows that, for all $x$,
\begin{equation*}
\sum_l \sum_{y\sim x} \left\vert c_{k,l}^j\right\vert \left\vert \chi_l^{j+1}(y)-\chi_l^{j+1}(x) \right\vert \leq C2^j.
\end{equation*}
Therefore,
$$
\sum_{x\in CB^j_k} \sum_{y\sim x} \left\vert \Delta_4(x,y)\right\vert^{q^{\ast}} m(x)
\leq C2^{(j+1)q^{\ast}} V(B_k^j).
$$
Gathering the estimates on $\Delta_i$, $1\leq i\leq 4$, we obtain \eqref{gradljk}.\par

\medskip

\noindent We now set $a^j_k= C^{-1} 2^{-j} V\left(B^j_k\right)^{-1} l^j_k$ and 
 $\lambda_{j,k}=C 2^j V\left(B^j_k\right).$ Then $f=\sum_{j,k}
  \lambda_{j,k}a^j_k,$ with $a^j_k$ being $(1,q^*)$ homogeneous Hardy-Sobolev
   atoms and 
   $$
   \begin{array}{lll}
  \displaystyle  \sum\limits_{j,k}\left\vert\lambda_{j,k}\right\vert&=&  \displaystyle C \sum\limits_{j,k} 2^j V\left(B^j_k\right)\\
   &\leq&  \displaystyle C \sum\limits_{j,k} 2^j V(\underline{B^j_k})\\
   &\leq&  \displaystyle C\sum\limits_j 2^j V\left(\{x: \mathcal{M}_q (Nf)(x) > 2^j\}\right)\\
   &\leq&  \displaystyle C \sum\limits_{x\in\Gamma} \mathcal{M}_q (Nf)(x) m(x)\\
   &\leq &  \displaystyle C_q \left\Vert Nf\right\Vert_{L^1(\Gamma)}\sim \left\Vert f\right\Vert_{\dot{S}^{1,1}(\Gamma)},
   \end{array}
   $$
where we used the fact that the $\underline{B^j_k}$ are pairwise disjoint.
\hfill\fin\par
\noindent{\bf Proof of Lemma \ref{observ}: } let $x\in B^j_k$ and $y\in B^{j+1}_l$ such that $x\sim y$. Denote by $x^j_k$ (resp. $x^{j+1}_l$) the center of $B^j_k$ (resp. $B^{j+1}_l$). Then
$$
d(x^j_k,x^{j+1}_l)\leq d(x^j_k,x)+d(x,y)+d(y,x^{j+1}_l)\leq r^j_k+r^{j+1}_l+1.
$$
Thus, since $F^j\subset F^{j+1}$,
$$
\begin{array}{lll}
\displaystyle r^{j+1}_l & = & \frac 12 d\left(x^{j+1}_l,F^{j+1}\right)\\
& \leq & \frac 12 d\left(x^{j+1}_l,x^j_k\right)+\frac 12 d(x^j_k,F^{j+1})\\
& \leq & \frac 12\left(r^j_k+r^{j+1}_l+1\right)+\frac 12 d(x^j_k,F^{j}),
\end{array}
$$
from which we deduce
$$
r^{j+1}_l\leq r^j_k+1+d(x^j_k,F^j)=r^j_k+1+2r^j_k\leq 4r^j_k,
$$
as claimed. The proof of $2$ is classical. \hfill\fin.
\subsection{Comparison between different atomic spaces}
In the present section, we show that $\dot{H}S^{1}_{t,ato}(\Gamma)=\dot{H}S^{1}_{t^{\prime},ato}(\Gamma)$ for all $t,t^{\prime}\in (1,+\infty]$, following ideas from \cite{badber}. 
We will need: 
\begin{lem}\label{the centered maximal function}
Assume that $\Gamma$ satisfies $(D)$.
\begin{itemize}
\item[$1.$]
Let 
\begin{equation*}
\mathcal{M}_c f(x) := \sup\limits_{r>0} \frac{1}{V(x,r)} \sum\limits_{B(x,r)} \vert f(y)\vert 
m(y)
\end{equation*}
be the centered maximal function of $f$. Observe that if $x\in B(y,r)$ then $B(y,r) 
\subset B(x,2r).$ It follows that 
\begin{equation*}
\mathcal{M}_c f \leq \mathcal{M}_{HL} f \leq C \mathcal{M}_c f
\end{equation*}
where $C$ only depends on the constant of the doubling property.
\item[$2.$]
Let $f$ be an $L^1$ function supported in $B_0=B(x_0,r_0).$ Then there is 
$C_1$ depending on the doubling constant such that 
\begin{equation*}
\Omega_{\alpha}:=\{x\in \Gamma: \mathcal{M}_{HL}(f)(x)> \alpha\}\subset B(x_0,2r_0)
\end{equation*}
\end{itemize}
whenever $\alpha > \frac{C_1}{V(B_0)} \sum\limits_{x\in B_0} \vert f(x) \vert m(x)$
\end{lem}
{\bf Proof: } it is obvious that ${\mathcal M}_cf\leq {\mathcal M}_{HL}f$ everywhere on $\Gamma$. Moreover, let $x\in \Gamma$ et $B=B(x_0,r)\ni x$ be a ball. Then $B\subset B(x,2r)\subset B(x_0,3r)$, so that
$$
\frac 1{V(B)} \sum_{y\in B} \left\vert f(y)\right\vert m(y)\leq \frac 1{V(x,2r)}\frac {V(x,2r)}{V(B)}  \sum_{y\in B(x,2r)} \left\vert f(y)\right\vert m(y)\leq C{\mathcal M}_cf(x),
$$
and the result follows by taking the supremum over all balls $B$ containing $x$.\par
\noindent For the second assertion, assume that $x\notin B(x_0,2r_0)$ and let $B=B(x,r)\ni x$ be a ball centered at $x$. Then
$$
\frac 1{V(B)}\sum_{y\in B} \left\vert f(y)\right\vert m(y)=\frac 1{V(B)} \sum_{y\in B\cap B(x_0,r_0)} \left\vert f(y)\right\vert m(y).
$$
If $B\cap B(x_0,r_0)=\emptyset$, then this quantity is $0$. Otherwise, $2r_0<d(x,x_0)\leq r+r_0$, so that $r_0<r$ and $B_0\subset B(x,2r)$. It follows that
$$
\frac 1{V(B)}\sum_{y\in B} \left\vert f(y)\right\vert m(y)\leq \frac 1{V(B_0)}\frac{V(B_0)}{V(B)}\sum_{y\in B_0} \left\vert f(y)\right\vert m(y)\leq \frac C{V(B_0)}\sum_{y\in B_0} \left\vert f(y)\right\vert m(y),
$$
which yields the conclusion by part $1.$, provided that $C_1$ is big enough. \hfill\fin\par
Let us now prove:
\begin{prop} \label{comparatom}
Let $\Gamma$ satisfying $(D)$ and the Poincar\'e inquality $(P_1)$. Then $HS^1_{t,ato}\subset HS^1_{\infty,ato}$ for every $t>1$ and 
therefore $HS^1_{t_1,ato}=HS^1_{t_2,ato}$ for every $1<t_1,t_2\leq +\infty$.
\end{prop}
{\bf Proof: } let $t>1$. It is enough to prove that there exists $C>0$ such that, for every $(1,t)-$atom $a$, $a$ belongs to $\dot{H}S^1_{\infty,ato}(\Gamma)$ with
$$
\left\Vert a\right\Vert_{\dot{H}S^1_{\infty,ato}(\Gamma)}\leq C.
$$
In the sequel, set ${\mathcal M}_{HL}^1:={\mathcal M}_{HL}$ and ${\mathcal M}_{HL}^{n+1}={\mathcal M}_{HL}^{n}\circ {\mathcal M}_{HL}$ for all $n\in \N$. Let $a$ be $(1,t)-$atom supported in a ball $B_0.$ Set $b=V(B_0)a.$
\newline
We claim that there exist $K,\alpha,C,N>0$ only depending on $t$ and the geometric constants with the following property: for all $l\in \N^{\ast}$, there exists a collection of balls $(B_{j_l})_{j_l\in \N^l}$ such 
that for every $n\geq 1$ 
\begin{equation}\label{definition of b}
b=CN \sum\limits_{l=1}^{n-1}(K\alpha)^{l+1}\sum\limits_{j_l\in \N^l} V(B_{j_l})a_{j_l}
+ \sum\limits_{j_n\in \N^n} h_{j_n}
\end{equation}
and, for all $n\in \N^{\ast}$,
\begin{equation} \label{prop1}
a_{j_l}\mbox{ is an }(1,\infty)\mbox{-atom supported in }B_{j_l}, 1\leq l\leq n-1,
\end{equation}
\begin{equation} \label{prop2}
\bigcup_{j_n\in \N^n} B_{j_n}\subset \Omega_l :=\left\{x\in \Gamma; \mathcal{M}_{HL}^{l+1}(\vert \nabla b \vert)(x)> K\frac{\alpha^l}{2}\right\},
\end{equation}
\begin{equation} \label{prop3}
\sum\limits_{j_l} {\bf 1}_{B_{j_l}} \leq N^l,
\end{equation}
\begin{equation} \label{prop4}
\mbox{ supp }h_{j_l} \subset B_{j_l}, \sum\limits_{x\in B_{j_l}} h_{j_l}(x)m(x)=0,
\end{equation}
\begin{equation} \label{prop5}
\left\vert \nabla h_{j_l}(x)\right\vert\leq C\left((\alpha K)^l \chi_{j_l}+\mathcal{M}_{HL}^l
(\vert \nabla 
b\vert)\right)(x)\mbox{ for all }x\in \Gamma,
\end{equation}
\begin{equation} \label{prop6}
\frac 1{V(B_{j_l})} \left\Vert \nabla h_{j_l}\right\Vert_{L^1(\Gamma)} \leq C(K\alpha)^l,
\end{equation}
where $\chi_{j_n}$ stands for the characteristic function of $B_{j_n}$. \par

\medskip

\noindent Let us assume that this construction is done. We claim that
\begin{equation}\label{decomposition of a}
a=\sum_{l=1}^{\infty} CN(K\alpha)^{l+1}\sum\limits_{j_l\in \N^l} \frac{V(B_{j_l})}{V(B_0)}a_{j_l},
\end{equation}
where the series converges in $\dot{W}^{1,1}(\Gamma)$ and
\begin{equation}\label{the bounded V(Bjn)}
\frac{N}{V(B_0)}\sum\limits_{l=1}^\infty (K\alpha)^{l+1} \sum\limits_{j_l\in \N^l} V
 (B_{j_l}) \leq C
\end{equation}
where $C$ is independent of $a$. \par
\noindent Let us first check \eqref{the bounded V(Bjn)}. Indeed, it follows from \eqref{prop1}, \eqref{prop3}  and the $L^t(\Gamma)$-boundedness of $\mathcal{M}_{HL}$ that 
\begin{equation*}
\sum\limits_{j_l} V(B_{j_l}) \leq C N^l m\left(\bigcup_{j_l}B_{j_l}\right) \leq C N^l m(\Omega_l)
\leq CN^l\left(\frac{2}{K \alpha^l}\right)^t \Vert \nabla b\Vert^t_{L^t(\Gamma)}.
\end{equation*}
As a consequence,
\begin{equation*}
\begin{array}{lll}
\displaystyle \sum_{l=0}^{\infty}(K\alpha)^l\sum_{j_l\in \N^l}V(B_{j_l})&\leq&\displaystyle C 2^t \sum_{n=0}^{\infty}(K\alpha)^l N^l (K\alpha^l)^{-t}\Vert \nabla b \Vert_{L^t(\Gamma)}^t\\
&\leq&\displaystyle C2^t K^{-t}\sum_{l=0}^{\infty}(N K \alpha^{(1-t)})^l \Vert \nabla b \Vert^t_{L^t(\Gamma)} 
\end{array}
\end{equation*}
and, since $\Vert \nabla b\Vert^t_{L^t(\Gamma)}\leq C V(B_0)$, we obtain \eqref{the bounded 
V(Bjn)} with $C$ only depending on $t,K,\alpha$ and $N$, provided that $\alpha$ is chosen such that 
$\frac{N K}{\alpha^{t-1}}<1.$\par
\noindent We now focus on \eqref{decomposition of a}. By \eqref{prop6}, one has
$$
\begin{array}{lll}
\displaystyle \frac 1{V(B_0)} \left\Vert \sum\limits_{j_n\in \N^n} h_{j_n}\right\Vert_{\dot{W}^{1,1}(\Gamma)} & \leq & \displaystyle \frac 1{V(B_0)} \sum_{j_n\in \N^n}Ê\left\Vert \nabla h_{j_n}\right\Vert_{L^1(\Gamma)}\\
& \leq & \displaystyle  C(K\alpha)^n \sum_{j_n\in \N^n}Ê\frac{V(B_{j_n})}{V(B_0)}
\end{array}
$$
and, by \eqref{the bounded 
V(Bjn)}, this quantity converges to $0$ when $n\rightarrow +\infty$, which yields \eqref{decomposition of a}. \par

\medskip

\noindent Let us now turn the the construction, which will be done by induction on $l$, starting with $l=1$. Set 
\begin{equation*}
\tilde{\Omega}_1=\{x\in \Gamma : \mathcal{M}_{HL}(\nabla b)(x) > K\alpha\},
\end{equation*}
where $K,\alpha$ will be chosen such that $K\alpha>C_1$ and $C_1$ is given by Lemma \ref{the centered maximal function}. Hence, $\tilde{\Omega}_1 \subset 
2B_0$. Moreover, 
$$
m\left(\widetilde{\Omega}_1\right) \leq \frac 1{(K\alpha)^t} \left\Vert \mathcal{M}_{HL}(\nabla b)\right\Vert_{L^t(\Gamma)}^t \leq \frac C{(K\alpha)^t} \left\Vert (\nabla b)\right\Vert_{L^t(\Gamma)}^t< +\infty.
$$
If $\widetilde{\Omega_1}=\emptyset$, then $\frac{b}
{NCK\alpha V(B_0)}$ is a $(1,\infty)$ atom and we are done. Assume now that $\widetilde{\Omega_1}\neq \emptyset$ and define the balls $B_i$ and the functions $\chi_i$ as in the proof of Proposition \ref{Calderon-Zygmund decomposition}. Set also
$$
h_{i}:=(b-c_i) \chi_{i},
$$
where
$$
c_i:=\frac 1{\sum_{x\in B_i} \chi_i(x)m(x)} \sum_{x\in B_i} b(x)\chi_i(x)m(x).
$$
Clearly, $\mbox{supp }h_i \subset B_i$. Moreover, 
\begin{equation} \label{cancelbi}
\sum_{x\in B_i} h_i(x)m(x)=0.
\end{equation}
We now claim:
\begin{equation} \label{gradhi}
\left\Vert \nabla h_i\right\Vert_{L^1(\Gamma)}\leq C\alpha V(B_i).
\end{equation}
Indeed, arguing as in the proof of Proposition \ref{decomposition f}, one has, for all $x\sim y\in \Gamma$,
$$
\begin{array}{lll}
\displaystyle b_i(y)-b_i(x) & = & \displaystyle ((b(y)-b(x))\chi_i(y)+ (b(x)-c_i)(\chi_i(y)-\chi_i(x))\\
& = & \displaystyle A(x,y)+B(x,y).
\end{array}
$$
On the one hand, using the support condition on $\chi_i$,
\begin{equation} \label{sumA}
\sum_{x\sim y} \left\vert A(x,y)\right\vert m(x) \leq \sum_{x\in 2B_i} \left\vert \nabla b(x)\right\vert m(x) \leq CV(B_i)K\alpha .
\end{equation}
On the other hand,
$$
\sum_{x\sim y} \left\vert B(x,y)\right\vert m(x) \leq \frac C{r_i} \sum_{x\in 2B_i} \left\vert b(x)-c_i\right\vert m(x).
$$
But
$$
\left\Vert b-c_i\right\Vert_{L^1(2B_i)}\leq \left\Vert b-b_{B_i}\right\Vert_{L^1(2B_i)} +CV(B_i) \left\vert b_{B_i}-c_i\right\vert,
$$
and, arguing as in the proof of Proposition \ref{decomposition f} and using $(P_1)$, one obtains
\begin{equation} \label{sumB}
\sum_{x\sim y} \left\vert B(x,y)\right\vert  m(x) \leq \sum_{x\in 2B_i} \left\vert \nabla b(x)\right\vert m(x) \leq CV(B_i)K\alpha.
\end{equation}
Thus, \eqref{sumA} and \eqref{sumB} yield \eqref{gradhi}.\par

\medskip

\noindent Define now the functions $g$ (denoted by $g_0$ in the sequel) and $h$ as in the proof of Proposition \ref{Calderon-Zygmund decomposition}, so that
\begin{equation}\label{decomposition b}
b=\sum\limits_j h_j + g_0.
\end{equation}
Observe that the series in \eqref{decomposition b} converges in $\dot{W}^{1,1}(\Gamma)$. Indeed, by \eqref{gradhi},
$$
\begin{array}{lll}
\displaystyle \sum_j \left\Vert \nabla h_j\right\Vert_{L^1(\Gamma)} & = & \displaystyle \sum_j \left\Vert \nabla h_j\right\Vert_{L^1(2B_j)}\\
& \leq & \displaystyle C\sum_j V(B_j)^{1-1/t}  \left\Vert \nabla h_j\right\Vert_{L^t(2B_j)}\\
& \leq & \displaystyle CK\alpha \sum_j V(B_j)\\
& \leq & \displaystyle C(K\alpha)^{1-t} \left\Vert \nabla b\right\Vert_{L^t(\Gamma)}^t\\
& \leq & \displaystyle C(K\alpha)^{1-t} V(B_0).
\end{array}
$$
Moreover, since $\sum b(x)m(x)=0$ and $\sum h_{j} (x)m(x)=0$ for all $j$, one also has $\sum g_0(x)m(x)=0$. \par
\noindent Arguing as in the proof of Proposition \ref{Calderon-Zygmund decomposition}, one establishes that 
$$
\left\Vert \nabla g_0\right\Vert_{L^{\infty}(\Gamma)}\leq CK\alpha.
$$
It follows that $a_0=\frac{g_0}
{NCK\alpha V(B_0)}$ 
is a $(1,\infty)$-atom, and \eqref{decomposition b} yields
\begin{equation*}
b=NCK\alpha V(B_0)a_0+\sum\limits_{j\in \N}h_j
\end{equation*}  
Thus, properties \eqref{prop1}, \eqref{prop2}, \eqref{prop3} and \eqref{prop4} hold. Property \eqref{prop6} has already been checked in \eqref{gradhi}. Moreover,
\begin{equation*}
\begin{array}{lll}
\vert \nabla h_j(x)\vert&\leq& \vert b(x)-c_j\vert \vert\nabla\chi_j(x)\vert +(\max
\limits_{y\sim x}\chi_j(y))\vert\nabla b(x)\vert\\
&=&I+II.
\end{array}
\end{equation*}
We estimate $I$ as follows:
$$
I\leq \frac C{r_j} \vert b(x)-c_j\vert.
$$
But, following the proof of Theorem 0.1 in \cite{badber}, if $l_j\in \Z$ is such that $2^{l_j}\leq r_j<2^{l_j+1}$, one has, using $(P_1)$, 
$$
\begin{array}{lll}
\displaystyle \left\vert b(x)-c_j\right\vert & \leq & \displaystyle \sum_{k=-(l_j+1)}^{-1} \left\vert b_{B(x,2^kr_j)}-b_{B(x,2^{k+1}r_j)}\right\vert + \left\vert b_{B(x,r_j)}-c_j\right\vert\\
& \leq & \displaystyle \sum_{k=-(l_j+1)}^{-1} \frac 1{V(x,2^kr_j)} \sum_{z\in B(x,2^kr_j)} \left\vert b(z)-b_{B(x,2^{k+1}r_j)} \right\vert m(z)+\left\vert b_{B(x,r_j)}-b_{2B_j}\right\vert\\
& + & \displaystyle \left\vert \frac 1{\sum_{z\in B_j} \chi_j(z)m(z)} \sum_{z\in B_j} \left(b(z)-\frac 1{V(2B_j)}\sum_{w\in 2B_j} b(w)m(w)\right)\chi_j(z)m(z)\right\vert \\
& \leq & \displaystyle  C \sum_{k=-(l_j+1)}^{-1} \frac 1{V(x,2^{k+1}r_j)} \sum_{z\in B(x,2^{k+1}r_j)} \left\vert b(z)-b_{B(x,2^{k+1}r_j)} \right\vert m(z)+\left\vert b_{B(x,r_j)}-b_{2B_j}\right\vert\\
& + & \displaystyle \left\vert \frac 1{\sum_{z\in B_j} \chi_j(z)m(z)} \sum_{z\in B_j} \left(b(z)-\frac 1{V(2B_j)}\sum_{w\in 2B_j} b(w)m(w)\right)\chi_j(z)m(z)\right\vert \\
& \leq & \displaystyle C \sum_{k=-(l_j+1)}^{-1} 2^{k+1}r_j {\mathcal M}_{HL}(\left\vert \nabla b\right\vert)(x)+ \frac 1{V(2B_j)} \sum_{z\in 2B_j} \left\vert b(z)-b_{2B_j}\right\vert m(z) \\
&+ & \displaystyle \frac 1{\sum_{z\in B_j} \chi_j(z)m(z)}\sum_{z\in 2B_j} \left\vert b(z)-\frac 1{V(2B_j)}\sum_{w\in 2B_j} b(w)m(w)\right\vert \left\vert \chi_j(z)\right\vert m(z)\\
& \leq & \displaystyle Cr_j\left( {\mathcal M}_{HL}(\left\vert \nabla b\right\vert)(x)+K\alpha\right).
\end{array}
$$
Moreover, $II\leq \vert \nabla b(x)\vert \leq \mathcal{M}_{HL}(\vert\nabla b \vert)(x)$. Finally, \eqref{prop5} is satisfied. The construction for $l=1$ is therefore complete. \par

\medskip

\noindent Assuming now that the construction is done for $l$, the construction for $l+1$ is performed by arguments analogous to the previous one (see also the proof of Theorem 0.1 in \cite{badber}). This ends the proof of Proposition \ref{comparatom}. \hfill\fin 
\subsection{Interpolation between Hardy-Sobolev and Sobolev spaces}
To establish Theorem \ref{interpol}, observe that, by Theorems \ref{characsobolev} and \ref{equal1}, $f\in \dot{S}^{1,1}(\Gamma)$ (resp. $f\in \dot{W}^{1,p}(\Gamma)$ if $p>1$) if and only if ${\mathcal M}^+f\in L^1(\Gamma)$ (resp. ${\mathcal M}^+f\in L^p(\Gamma)$. Therefore, Theorem \ref{interpol} follows from the classical 
linearization method of maximal operators (see \cite{steinweiss}, Chapter 5).
\section{Boundedness of Riesz transforms} \label{Rieszproofs}
\subsection{The boundedness of Riesz transforms on Hardy-Sobolev spaces}
This section is devoted to the proof of Theorem \ref{Riesz}. We first establish:
\begin{prop} \label{propestimatom}
There exists $C>0$ such that, for all atom $a\in H^1(\Gamma)$, $(I-P)^{-1/2}a\in \dot{S}^{1,1}(\Gamma)$ and
\begin{equation} \label{estimatom}
\left\Vert (I-P)^{-1/2}a\right\Vert_{\dot{S}^{1,1}(\Gamma)}\leq C.
\end{equation}
\end{prop}
The proof relies on some estimates for the iterates of $p$, taken from \cite{scand,pota}. Define
\[
p_0(x,y):=
\left\{
\begin{array}{ll}
1 &\mbox{ if }x=y,\\
0 & \mbox{ if }x\neq y,
\end{array}
\right.
\]
and, for all $k\in \N$ and all $x,y\in \Gamma$,
\[
p_{k+1}(x,y)=\sum\limits_{z\in \Gamma}p(x,z)p_k(z,y).
\]
By \eqref{reversibility}, one has
$$
p_k(x,y)m(x)=p_k(y,x) m(y)
$$
for all $k\in \N$ and all $x,y\in \Gamma$. \par
\noindent Let $y_0\in \Gamma$. For all $k\in \N$ and all $x\in \Gamma$, define
$$
q_k(x,y):=\frac{p_k(y,x)-p_k(y_0,x)}{m(x)}.
$$
Recall the following bounds on $p_k$ and $q_k$ (\cite{scand}, Lemmata 2 and 4 and \cite{pota}, Lemmata 28 and 29):
\begin{lem} \label{estimpk}
There exist $C,\alpha>0$ such that, for all $y\in \Gamma$,
\begin{itemize}
\item[$1.$]
$$
\sum_{x\in \Gamma} \left\vert \nabla_x p_k(x,y)\right\vert^2 \exp\left(\alpha\frac{d^2(x,y)}k\right)m(x)\leq \frac C{V(y,\sqrt{k})}m^2(y),
$$
\item[$2.$]
$$
\sum_{x\in \Gamma} \left\vert \nabla_x p_k(x,y)\right\vert^2 \exp\left(\alpha\frac{d^2(x,y)}k\right)m(x)\leq \frac C{kV(y,\sqrt{k})}m^2(y).
$$
\end{itemize}
\end{lem} 
\begin{lem} \label{estimqk}
There exist $C,h,\alpha>0$ such that, for all $y_0,y\in \Gamma$ and all $k\geq 1$ such that $d(y,y_0)\leq \sqrt{k}$,
\begin{itemize}
\item[$1.$]
$$
\sum_{x\in \Gamma} \left\vert q_k(x,y)\right\vert^2\exp\left(\alpha\frac{d^2(x,y)}k\right)m(x)\leq \frac C{V(y,\sqrt{k})}\left(\frac{d(y,y_0)}{\sqrt{k}}\right)^h,
$$
\item[$2.$]
$$
\sum_{x\in \Gamma} \left\vert \nabla_xq_k(x,y)\right\vert^2\exp\left(\alpha\frac{d^2(x,y)}k\right)m(x)\leq \frac C{kV(y,\sqrt{k})}\left(\frac{d(y,y_0)}{\sqrt{k}}\right)^h.
$$
\end{itemize}
\end{lem}
{\bf Proof of Proposition \ref{propestimatom}: } let $a$ be an atom supported in $B=B(y_0,r)$. Pick up a sequence of functions $(\chi_j)_{j\geq 0}$ such that 
$$
\mbox{supp }\chi_0\in 4B,\mbox{ supp }\chi_j\subset 2^{j+2}B\setminus 2^{j-1}B,\ \left\Vert d\chi_j\right\Vert_{\infty}\leq \frac C{2^jr}
$$
and
$$
\sum_{j\geq 0} \chi_j=1\mbox{ on }\Gamma.
$$
For all $j\geq 0$, all $x\in \Gamma$ and all $y\sim x$, one has
$$
\begin{array}{lll}
\displaystyle \chi_j(y)(I-P)^{-1/2}a(y)-\chi_j(x)(I-P)^{-1/2}a(x) & = & \displaystyle \chi_j(y)\left((I-P)^{-1/2}a(y)-(I-P)^{-1/2}a(x)\right)\\
& + & \displaystyle (I-P)^{-1/2}a(x)\left(\chi_j(y)-\chi_j(x)\right).
\end{array}
$$
It follows that, if $\nabla \left(\chi_j(I-P)^{-1/2}a\right)(x)\neq 0$, then either $\chi_j(x)\neq 0$, or there exists $y\sim x$ such that $\chi_j(y)\neq 0$. As a consequence, $\mbox{ supp } \nabla \left(\chi_j(I-P)^{-1/2}a\right)\subset C_j(B):=2^{j+3}B\setminus 2^{j-2}B$ if $j\geq 3$ and $\mbox{ supp } \nabla \left(\chi_j(I-P)^{-1/2}a\right)\subset C_j(B):=2^{j+3}B$ if $j\leq 2$. Decompose $(I-P)^{-1/2}a$ as
$$
\begin{array}{lll}
\displaystyle (I-P)^{-1/2}a & = & \displaystyle \sum_{j\geq 0} \chi_j(I-P)^{-1/2}a\\
& = & \displaystyle \sum_{j\geq 0} V^{1/2}(2^{j+3}B)\left\Vert \nabla \left(\chi_j(I-P)^{-1/2}a\right)\right\Vert_{L^2(\Gamma)} \frac{\chi_j(I-P)^{-1/2}a}{V^{1/2}(2^{j+3}B)\left\Vert \nabla \left(\chi_j(I-P)^{-1/2}a\right)\right\Vert_{L^2(\Gamma)}}\\
& := & \displaystyle \sum_{j\geq 0} V^{1/2}(2^{j+3}B)\left\Vert \nabla \left(\chi_j(I-P)^{-1/2}a\right)\right\Vert_{L^2(\Gamma)} b_j.
\end{array}
$$
We first check that, for all $j\geq 0$, up to a constant only depending on the constants of the graph $\Gamma$, $b_j$ is an atom in $\dot{H}S^{1}_{2,ato}(\Gamma)$ if, in Definition \ref{Atom}, condition $3$ is replaced by condition $3^{\prime}$ in Remark \ref{relax}. Indeed, since \eqref{doubling} and $(P_1)$ hold, there exists $C>0$ such that, for all balls $B$ of radius $r$ and all functions $f\in W^{1,2}_0(B)$,
\begin{equation} \label{poincaredirichlet}
\left\Vert f\right\Vert_{L^2(B)}\leq Cr\left\Vert \nabla f\right\Vert_{L^2(B)}
\end{equation}
(see \cite{BR09}, inequality $(8.2)$). Then, for all $j$, since $\chi_j$ is supported in $2^{j+2}B$, \eqref{poincaredirichlet} yields
$$
\left\Vert \chi_j(I-P)^{-1/2}a\right\Vert_{L^2(2^{j+2}B)}\leq C2^{j+2}r \left\Vert \nabla \left(\chi_j(I-P)^{-1/2}a\right)\right\Vert_{L^2(\Gamma)},
$$
which shows that
$$
\left\Vert b_j\right\Vert_{L^2(2^{j+2}B)}\leq C2^{j+2}rV^{-1/2}(2^{j+2}B),
$$
as claimed.\par
The estimate \eqref{estimatom} will therefore be a consequence of 
\begin{equation} \label{toprove}
\sum_{j\geq 0} V^{1/2}(2^{j+3}B)\left\Vert \nabla \left(\chi_j(I-P)^{-1/2}a\right)\right\Vert_{L^2(\Gamma)}\leq C.
\end{equation}
Write
$$
\begin{array}{lll}
\displaystyle \left\Vert \nabla \left(\chi_j(I-P)^{-1/2}a\right)\right\Vert_{L^2(\Gamma)} & \leq & \left\Vert \nabla (I-P)^{-1/2}a\right\Vert_{L^2(C_j(B))}+\left\Vert \nabla \chi_j\right\Vert_{\infty} \left\Vert (I-P)^{-1/2}a\right\Vert_{L^2(C_j(B))}\\
& := & S_j+T_j.
\end{array}
$$
Let us first focus on $T_j$. As in \cite{BR09}, we use the expansion
$$
\begin{array}{lll}
\displaystyle (I-P)^{-1/2}a & = & \displaystyle \sum_{k=0}^{+\infty} a_kP^ka\\
& = & \displaystyle  \sum_{k=0}^{r^2} a_kP^ka +  \sum_{k=r^2+1}^{+\infty} a_kP^ka\\
& := & f_1+f_2,
\end{array}
$$
where the $a_k$'s are defined by
$$
(1-x)^{-1/2}=\sum_{k=0}^{+\infty} a_kx^k
$$
for all $x\in (-1,1)$. Recall that, when $k\rightarrow +\infty$,
\begin{equation} \label{equivak}
a_k\sim \frac 1{\sqrt{k\pi}}.
\end{equation}
For $f_1$, 
\begin{equation} \label{g1l2}
\left\Vert f_1\right\Vert_{L^2(C_j(B))}\leq \sum_{k=0}^{r^2} a_k \left\Vert P^ka\right\Vert_{L^2(C_j(B))}.
\end{equation}
For $k=0$, $P^ka=a$ so that
$$
\left\Vert P^ka\right\Vert_{L^2(\Gamma)}\leq V(B)^{-1/2}.
$$
Let $h\in L^2(C_j(B))$ with $\left\Vert h\right\Vert_{L^2}\leq 1$. For all $1\leq k\leq r^2$, Lemma \ref{estimpk} yields
\begin{equation} \label{estimh}
\begin{array}{lll}
\displaystyle \left\vert \sum_{x\in C_j(B)} P^ka(x) h(x)m(x)\right\vert & \leq  & \displaystyle \sum_{x\in C_j(B)} \left\vert h(x)\right\vert  \left(\sum_{y\in \Gamma} p_k(x,y)\left\vert a(y)\right\vert\right)m(x)\\
& = & \displaystyle \sum_{y\in \Gamma}Ê\left\vert a(y) \right\vert \left(\sum_{x\in C_j(B)} p_k(x,y) \exp\left(\frac{\alpha d^2(x,y)}{2k}\right)\exp\left(-\frac{\alpha d^2(x,y)}{2k}\right)\right.\\
& & \left.  \left\vert h(x)\right\vert  m(x)\right)\\
& \leq & \displaystyle e^{-c\frac{2^{2j}r^2}{k}}\sum_{y\in \Gamma}Ê\left\vert a(y) \right\vert \left(\sum_{x\in 2^{j+3}B} \left\vert p_k(x,y)\right\vert^2 \exp\left(\frac{\alpha d^2(x,y)}{k}\right)m(x)\right)^{1/2} \\
& & \displaystyle \left\Vert h\right\Vert_{L^2(C_j(B))}\\
& \leq & \displaystyle C e^{-c\frac{2^{2j}r^2}{k}} \sum_{y\in B}Ê\frac{\left\vert a(y) \right\vert}{V^{1/2}(y,\sqrt{k})} m(y).
\end{array}
\end{equation}
But, for all $y\in B$, \eqref{doubling} shows that
\begin{equation}Ê\label{1survol}
\begin{array}{lll}
\displaystyle \frac 1{V(y,\sqrt{k})} & = & \displaystyle \frac 1{V(y_0,\sqrt{k})} \frac {V(y_0,\sqrt{k})}{V(y,\sqrt{k})}\\
& \leq & \displaystyle \frac 1{V(y_0,\sqrt{k})} \frac {V(y,\sqrt{k}+r)}{V(y,\sqrt{k})}\\
& \leq & \displaystyle \frac 1{V(y_0,\sqrt{k})} \left(1+\frac r{\sqrt{k}}\right)^D\\
& \leq & \displaystyle \frac 1{V(2^{j+3}B)} \frac{V(2^{j+3}B)}{V(y_0,\sqrt{k})} \left(1+\frac r{\sqrt{k}}\right)^D\\
& \leq & \displaystyle \frac 1{V(2^{j+3}B)} \left(1+\frac {2^{j+3}r}{\sqrt{k}}\right)^{2D}.
\end{array}
\end{equation}
Therefore, it follows from \eqref{estimh} and the fact that $\left\Vert a\right\Vert_1\leq 1$ that
$$
\displaystyle \left\Vert P^ka \right\Vert_{L^2(C_j(B))} \leq  \frac C{V^{1/2}(2^{j+3}B)}  \exp\left(-c^{\prime}\frac{2^{2j}r^2}{k}\right). 
$$
Since, when $j\geq 3$ and $k=0$, $P^ka=a$ and $C_j(B)$ are disjoint, one obtains
\begin{equation} \label{estimg1}
\left\Vert f_1\right\Vert_{L^2(C_j(B))} \leq \frac C{V^{1/2}(2^{j+3}B)}\left(\sum_{k=1}^{r^2} \frac 1{\sqrt{k}} \exp\left(c^{\prime}\frac{2^{2j}r^2}{k}\right)+c_j\right),
\end{equation}
with $c_j=1$ if $j\leq 2$ and $c_j=0$ if $j\geq 3$. \par
\noindent For $f_2$, 
$$
\left\Vert f_2\right\Vert_{L^2(C_j(B))}\leq \sum_{k=r^2+1}^{\infty} a_k \left\Vert \sum_{y\in \Gamma} q_k(\cdot,y)a(y)m(y)\right\Vert_{L^2(C_j(B))}.
$$
Pick up a function $h\in L^2(C_j(B))$ with $\left\Vert h\right\Vert_{L^2}\leq 1$ again. For all $k\geq r^2+1$, Lemma \ref{estimqk} yields
\begin{equation} \label{estimhbis}
\begin{array}{lll}
\displaystyle \left\vert \sum_{x\in C_j(B)} P^ka(x) h(x)m(x)\right\vert & \leq  & \displaystyle \sum_{x\in C_j(B)} \left\vert h(x)\right\vert  \left(\sum_{y\in \Gamma} q_k(x,y)\left\vert a(y)\right\vert m(y)\right)m(x)\\
& = & \displaystyle \sum_{y\in \Gamma}Ê\left\vert a(y) \right\vert \left(\sum_{x\in C_j(B)} q_k(x,y) \exp\left(\frac{\alpha d^2(x,y)}{2k}\right)\exp\left(-\frac{\alpha d^2(x,y)}{2k}\right)\right.\\
& & \left.  \left\vert h(x)\right\vert  m(x)\right)m(y)\\
& \leq & \displaystyle e^{-c\frac{2^{2j}r^2}{k}}\sum_{y\in \Gamma}Ê\left\vert a(y) \right\vert \left(\sum_{x\in C_j(B)} \left\vert q_k(x,y)\right\vert^2 \exp\left(\frac{\alpha d^2(x,y)}{k}\right)m(x)\right)^{1/2} \\
& \times & \displaystyle \left\Vert h\right\Vert_{L^2(C_j(B))}m(y)\\
& \leq & \displaystyle Ce^{-c\frac{2^{2j}r^2}{k}}\left(\frac r{\sqrt{k}}\right)^{h/2} \sum_{y\in B}Ê\frac{\left\vert a(y) \right\vert}{V^{1/2}(y,\sqrt{k})} m(y)\\
& \leq & \displaystyle \frac C{V^{1/2}(y_0,\sqrt{k})} e^{-c\frac{2^{2j}r^2}{k}} \left(\frac r{\sqrt{k}}\right)^{h/2}.
\end{array}
\end{equation}
Arguing as before and using \eqref{reversevol}, one therefore obtains
\begin{equation} \label{estimg2}
\begin{array}{lll}
\displaystyle \left\Vert f_2\right\Vert_{L^2(C_j(B))} & \leq  & \displaystyle \sum_{k=r^2+1}^{+\infty} \frac C{\sqrt{k}V^{1/2}(y_0,\sqrt{k})} e^{-c\frac{2^{2j}r^2}{k}} \left(\frac r{\sqrt{k}}\right)^{h/2}\\
& = & \displaystyle \frac C{V^{1/2}(2^{j+3}B)} \sum_{k=r^2+1}^{+\infty}  \frac 1{\sqrt{k}} \frac{V^{1/2}(2^{j+3}B)}{V^{1/2}(y_0,\sqrt{k})} e^{-c\frac{2^{2j}r^2}{k}} \left(\frac r{\sqrt{k}}\right)^{h/2}\\
& \leq & \displaystyle \frac C{V^{1/2}(2^{j+3}B)} \sum_{k=r^2+1}^{+\infty} \frac 1{\sqrt{k}} f\left(\frac{2^{j+3}r}{\sqrt{k}}\right) e^{-c\frac{2^{2j}r^2}{k}} \left(\frac r{\sqrt{k}}\right)^{h/2},
\end{array}
\end{equation}
where
$$
f(u)=\left\{
\begin{array}{ll}
u^{D/2} & \mbox{ if }u>1,\\
u^{d/2} & \mbox{ if }u\leq 1.
\end{array}
\right.
$$
Gathering \eqref{estimg1} and \eqref{estimg2}, one therefore obtains
$$
T_j \leq \frac C{2^jrV^{1/2}(2^{j+3}B)} \left(c_j+\sum_{k=1}^{r^2} \frac 1{\sqrt{k}} e^{-c\frac{2^{2j}r^2}{k}} + \sum_{k=r^2+1}^{+\infty} \frac 1{\sqrt{k}} e^{-c\frac{2^{2j}r^2}{k}} f\left(\frac{2^{j+3}r}{\sqrt{k}}\right) \left(\frac r{\sqrt{k}}\right)^{h/2} \right).
$$
Thus,
$$
\begin{array}{lll}
\displaystyle V^{1/2}(2^{j+3}B) T_j & \leq & \displaystyle \frac C{2^jr}\left(c_j+\sum_{k=1}^{r^2} \frac 1{\sqrt{k}} e^{-c\frac{2^{2j}r^2}{k}} + \sum_{k=r^2+1}^{+\infty} \frac 1{\sqrt{k}} e^{-c\frac{2^{2j}r^2}{k}} f\left(\frac{2^{j+3}r}{\sqrt{k}}\right) \left(\frac r{\sqrt{k}}\right)^{h/2}\right) \\
& \leq & \displaystyle \frac C{2^jr} \int_0^{2r^2}e^{-c\frac{2^{2j}r^2}{t}}\frac{dt}{\sqrt{t}}+  \int_{r^2}^{+\infty} e^{-c\frac{2^{2j}r^2}{t}}  f\left(\frac{2^{j+3}r}{\sqrt{t}}\right) \left(\frac r{\sqrt{t}}\right)^{h/2} \frac{dt}t\\
& = & \displaystyle C\int_{2^{2j-1}}^{+\infty} e^{-cu} \frac{du}{u^{3/2}}+C\int_0^{2^{2j}} e^{-cu} f(8\sqrt{u})\left(\frac{\sqrt{u}}{2^j}\right)^{h/2}\frac{du}{u^{3/2}}\\
& \leq & \displaystyle C2^{-jh/2}.
\end{array}
$$
Note that we used the fact that $d\geq 1$ in the last inequality (this is the only place where this assumption is used). Finally,
\begin{equation} \label{toprovetj}
\sum_{j=0}^{+\infty} V^{1/2}(2^{j+3}B) T_j \leq C.
\end{equation}

\noindent Let us now focus on $S_j$. For $j\leq 2$, the $L^2$-boundedness of $\nabla(I-P)^{-1/2}$ yields
$$
\left\Vert \nabla (I-P)^{-1/2}a\right\Vert_{L^2(C_j(B))}\leq \left\Vert \nabla (I-P)^{-1/2}a\right\Vert_{L^2(\Gamma)}\leq CV(B)^{-1/2}.
$$
Take now $j\geq 3$. As before, one has
$$
\begin{array}{lll}
\displaystyle \nabla (I-P)^{-1/2}a &\leq & \displaystyle \sum_{k=0}^{+\infty} a_k\nabla P^k a\\
& = & \displaystyle \sum_{k=0}^{r^2} a_k\nabla P^k a + \sum_{k=r^2+1}^{+\infty}Êa_k\nabla P^k a\\
& := & \displaystyle g_1+g_2.
\end{array}
$$
We estimate the $L^2$-norms of $g_1$ and $g_2$. For $g_1$, 
\begin{equation} \label{f1l2}
\left\Vert g_1\right\Vert_{L^2(C_j(B))}\leq \sum_{k=0}^{r^2} a_k \left\Vert \nabla P^ka\right\Vert_{L^2(C_j(B))}.
\end{equation}
Notice that, when $k=0$, $\nabla P^ka=\nabla a$ is supported in $2B$, which is disjoint from $C_j(B)$ since $j\geq 3$. Let $h\in L^2(C_j(B))$ with $\left\Vert h\right\Vert_{L^2}\leq 1$. For all $1\leq k\leq r^2$, Lemma \ref{estimpk} yields
\begin{equation} \label{estimg}
\begin{array}{lll}
\displaystyle \left\vert \sum_{x\in C_j(B)} \nabla P^ka(x) g(x)m(x)\right\vert & \leq  & \displaystyle \sum_{x\in C_j(B)} \left\vert h(x)\right\vert  \left(\sum_{y\in \Gamma} \nabla_xp_k(x,y)\left\vert a(y)\right\vert\right)m(x)\\
& = & \displaystyle \sum_{y\in \Gamma}Ê\left\vert a(y) \right\vert \left(\sum_{x\in C_j(B)} \nabla_xp_k(x,y) \exp\left(\frac{\alpha d^2(x,y)}{2k}\right)\exp\left(-\frac{\alpha d^2(x,y)}{2k}\right)\right.\\
& & \left.  \left\vert h(x)\right\vert  m(x)\right)\\
& \leq & \displaystyle e^{-c\frac{2^{2j}r^2}{k}}\sum_{y\in \Gamma}Ê\left\vert a(y) \right\vert \left(\sum_{x\in C_j(B)} \left\vert \nabla_xp_k(x,y)\right\vert^2 \exp\left(\frac{\alpha d^2(x,y)}{k}\right)m(x)\right)^{1/2} \\
& & \displaystyle \left\Vert h\right\Vert_{L^2(C_j(B))}\\
& \leq & \displaystyle \frac C{\sqrt{k}}e^{-c\frac{2^{2j}r^2}{k}} \sum_{y\in B}Ê\frac{\left\vert a(y) \right\vert}{V^{1/2}(y,\sqrt{k})} m(y).
\end{array}
\end{equation}
Thus, it follows from \eqref{1survol}, \eqref{estimg} and the fact that $\left\Vert a\right\Vert_1\leq 1$ that
$$
\displaystyle \left\Vert \nabla P^ka \right\Vert_{L^2(C_j(B))} \leq  \frac C{\sqrt{k}V^{1/2}(2^{j+3}B)}  \exp\left(-c^{\prime}\frac{2^{2j}r^2}{k}\right).
$$
As a consequence of \eqref{equivak} and \eqref{f1l2}, one therefore has
\begin{equation} \label{estimf1}
\left\Vert g_1\right\Vert_{L^2(C_j(B))} \leq \frac C{V^{1/2}(2^{j+3}B)}\left(c_j+\sum_{k=1}^{r^2} \frac 1{k}  \exp\left(c^{\prime}\frac{2^{2j}r^2}{k}\right)\right),
\end{equation}
where, again, $c_j=1$ if $j\leq 2$ and $c_j=0$ if $j\geq 3$.\par
\noindent For $g_2$, observe that, for all $x\in\Gamma$, since $\sum_{y\in \Gamma} a(y)m(y)=0$,
$$
\begin{array}{lll}
\displaystyle P^ka(x) & = & \displaystyle \sum_{y\in \Gamma} \frac{p_k(x,y)}{m(y)}a(y)m(y)\\
& = & \displaystyle \sum_{y\in \Gamma} \left(\frac{p_k(x,y)}{m(y)}-\frac{p_k(x,y_0)}{m(y_0)}\right)a(y)m(y)\\
& = & \displaystyle \frac 1{m(x)}\sum_{y\in \Gamma} \left(p_k(y,x)-p_k(y_0,x)\right)a(y)m(y)\\
& = & \displaystyle \sum_{y\in \Gamma}Êq_k(x,y)a(y)m(y).
\end{array}
$$
As a consequence,
$$
\left\Vert g_2\right\Vert_{L^2(C_j(B))}\leq \sum_{k=r^2+1}^{\infty} a_k \left\Vert \sum_{y\in \Gamma} \nabla_xq_k(\cdot,y)a(y)m(y)\right\Vert_{L^2(2^{j+3}B)}.
$$
Pick up a function $h\in L^2(C_j(B))$ with $\left\Vert h\right\Vert_{L^2}\leq 1$ again. For all $k\geq r^2+1$, Lemma \ref{estimqk} yields
\begin{equation} \label{estimgbis}
\begin{array}{lll}
\displaystyle \left\vert \sum_{x\in C_j(B)} \nabla P^ka(x) h(x)m(x)\right\vert & \leq  & \displaystyle \sum_{x\in C_j(B)} \left\vert h(x)\right\vert  \left(\sum_{y\in \Gamma} \nabla_xq_k(x,y)\left\vert a(y)\right\vert m(y)\right)m(x)\\
& = & \displaystyle \sum_{y\in \Gamma}Ê\left\vert a(y) \right\vert \left(\sum_{x\in C_j(B)} \nabla_xq_k(x,y) \exp\left(\frac{\alpha d^2(x,y)}{2k}\right)\exp\left(-\frac{\alpha d^2(x,y)}{2k}\right)\right.\\
& & \left.  \left\vert h(x)\right\vert  m(x)\right)m(y)\\
& \leq & \displaystyle e^{-c\frac{2^{2j}r^2}{k}}\sum_{y\in \Gamma}Ê\left\vert a(y) \right\vert \left(\sum_{x\in C_j(B)} \left\vert \nabla_xq_k(x,y)\right\vert^2 \exp\left(\frac{\alpha d^2(x,y)}{k}\right)m(x)\right)^{1/2} \\
& \times & \displaystyle \left\Vert h\right\Vert_{L^2(C_j(B))}m(y)\\
& \leq & \displaystyle \frac C{\sqrt{k}}e^{-c\frac{2^{2j}r^2}{k}}\left(\frac r{\sqrt{k}}\right)^{h/2} \sum_{y\in B}Ê\frac{\left\vert a(y) \right\vert}{V^{1/2}(y,\sqrt{k})} m(y)\\
& \leq & \displaystyle \frac C{\sqrt{k}V^{1/2}(y_0,\sqrt{k})} e^{-c\frac{2^{2j}r^2}{k}} \left(\frac r{\sqrt{k}}\right)^{h/2}.
\end{array}
\end{equation}
Using \eqref{equivak} again, as well as \eqref{doubling} and \eqref{reversevol}, one therefore obtains
\begin{equation} \label{estimf2}
\begin{array}{lll}
\displaystyle \left\Vert g_2\right\Vert_{L^2(C_j(B))} & \leq  & \displaystyle \sum_{k=r^2+1}^{+\infty} \frac C{kV^{1/2}(y_0,\sqrt{k})} e^{-c\frac{2^{2j}r^2}{k}} \left(\frac r{\sqrt{k}}\right)^{h/2}\\
& = & \displaystyle \frac C{V^{1/2}(2^{j+3}B)} \sum_{k=r^2+1}^{+\infty}  \frac 1k \frac{V^{1/2}(2^{j+3}B)}{V^{1/2}(y_0,\sqrt{k})} e^{-c\frac{2^{2j}r^2}{k}} \left(\frac r{\sqrt{k}}\right)^{h/2}\\
& \leq & \displaystyle \frac C{V^{1/2}(2^{j+3}B)} \sum_{k=r^2+1}^{+\infty} \frac 1k f\left(\frac{2^{j+3}r}{\sqrt{k}}\right) e^{-c\frac{2^{2j}r^2}{k}} \left(\frac r{\sqrt{k}}\right)^{h/2},
\end{array}
\end{equation}
Gathering \eqref{estimf1} and \eqref{estimf2}, one therefore obtains
$$
S_j \leq \frac C{V^{1/2}(2^{j+3}B)} \left(c_j+\sum_{k=1}^{r^2} \frac 1k e^{-c\frac{2^{2j}r^2}{k}} + \sum_{k=r^2+1}^{+\infty} \frac 1k e^{-c\frac{2^{2j}r^2}{k}} f\left(\frac{2^{j+3}r}{\sqrt{k}}\right) \left(\frac r{\sqrt{k}}\right)^{h/2} \right).
$$
Thus,
$$
\begin{array}{lll}
\displaystyle V^{1/2}(2^{j+3}B) S_j & \leq & \displaystyle C\left(c_j+\sum_{k=1}^{r^2} \frac 1k e^{-c\frac{2^{2j}r^2}{k}} + \sum_{k=r^2+1}^{+\infty} \frac 1k e^{-c\frac{2^{2j}r^2}{k}} f\left(\frac{2^{j+3}r}{\sqrt{k}}\right) \left(\frac r{\sqrt{k}}\right)^{h/2}\right) \\
& \leq & \displaystyle C\int_0^{2r^2}e^{-c\frac{2^{2j}r^2}{t}}\frac{dt}t+ C\int_{r^2}^{+\infty} e^{-c\frac{2^{2j}r^2}{t}}  f\left(\frac{2^{j+3}r}{\sqrt{t}}\right) \left(\frac r{\sqrt{k}}\right)^{h/2} \frac{dt}t\\
& = & \displaystyle C\int_{2^{2j-1}}^{+\infty} e^{-cu}\frac{du}u+C\int_0^{2^{2j}} e^{-cu} f(8\sqrt{u})\left(\frac{\sqrt{u}}{2^j}\right)^{h/2}\frac{du}u\\
& \leq & \displaystyle C2^{-jh/2},
\end{array}
$$
which proves that
\begin{equation} \label{toprovesj}
\sum_{j=0}^{+\infty} V^{1/2}(2^{j+3}B) S_j \leq C.
\end{equation}
Finally, \eqref{toprovesj} and \eqref{toprovetj} yield \eqref{toprove} and the proof of Proposition \ref{propestimatom} is complete. \hfill\fin\par
\noindent Let us now derive Theorem \ref{Riesz} from Proposition \ref{propestimatom}. Take $f\in H^1(\Gamma)$ and decompose
$$
f=\sum_{j=0}^{+\infty} \lambda_ja_j
$$
with $\sum_{j=0}^{+\infty} \left\vert \lambda_j\right\vert\leq 2\left\Vert f\right\Vert_{H^1(\Gamma)}$. For all $J\geq 0$, define
$$
f_J:=\sum_{j=0}^J \lambda_ja_j,
$$
so that $f_J\rightarrow f$ in $H^1(\Gamma)$. For all $j_1<j_2$,
$$
(I-P)^{-1/2}f_{j_2}-(I-P)^{-1/2}f_{j_1}=\sum_{j_1<j\leq j_2} \lambda_j (I-P)^{-1/2}a_j,
$$
which entails, by Proposition \ref{estimatom},
$$
\begin{array}{lll}
\displaystyle \left\Vert (I-P)^{-1/2}f_{j_2}-(I-P)^{-1/2}f_{j_1}\right\Vert_{\dot{S}^{1,1}(\Gamma)} & \leq  & \displaystyle \sum_{j_1<j\leq j_2} \left\vert \lambda_j\right\vert  \left\Vert (I-P)^{-1/2}a_j\right\Vert_{\dot{S}^{1,1}(\Gamma)}\\
& \leq & \displaystyle C\sum_{j_1<j\leq j_2} \left\vert \lambda_j\right\vert.
\end{array}
$$
This shows that $((I-P)^{-1/2}f_j)_{j\geq 0}$ is a Cauchy sequence in $\dot{S}^{1,1}(\Gamma)$, and therefore converges to some function $g\in \dot{S}^{1,1}(\Gamma)$. Moreover, using Proposition \ref{propestimatom} again,
$$
\left\Vert g\right\Vert_{\dot{S}^{1,1}(\Gamma)}=\lim_{J\rightarrow +\infty} \left\Vert (I-P)^{-1/2}f_J\right\Vert_{\dot{S}^{1,1}(\Gamma)}\leq C\sum_{j=0}^J \left\vert \lambda_j\right\vert\leq 2C\left\Vert f\right\Vert_{H^1(\Gamma)}.
$$
Furthermore, since $f_J\rightarrow f$ in $H^1(\Gamma)$, $d(I-P)^{-1/2}f_J\rightarrow d(I-P)^{-1/2}f$ in $L^{1}(E)$ (see \cite{pota}, Theorem 2.1). Since $d(I-P)^{-1/2}f_J\rightarrow dg$ in $L^1(E)$ by what we have just proved, $d(I-P)^{-1/2}f=dg$. As a consequence, $g=(I-P)^{-1/2}f\in \dot{S}^{1,1}(\Gamma)$ and
$$
\left\Vert (I-P)^{-1/2}f\right\Vert_{\dot{S}^{1,1}(\Gamma)}\leq 2C\left\Vert f\right\Vert_{H^1(\Gamma)},
$$
which concludes the proof of Theorem \ref{Riesz}. \hfill\fin
\subsection{Riesz transforms and Hardy spaces on edges}
Apart from Theorem \ref{Riesz}, it is also possible to establish that the Riesz transform maps $H^1(\Gamma)$ into a Hardy space on $E$, under assumptions \eqref{doubling} and $(P_1)$, without assuming \eqref{reversevol}. \par
\noindent Indeed, since $E$, endowed with its distance $d$ and its measure $\mu$, is also a space of homogeneous type (see Section \ref{measures}), we can define an atomic Hardy space on $E$. More precisely, an atom is a function $A\in L^2(E,\mu)$ (recall that $A$ is antisymmetric), supported in a ball $B\subset E$ and satisfying
$$
\sum_{(x,y)\in B} A(x,y)\mu_{xy}=0 \mbox{ and } \left\Vert A\right\Vert_{L^2(E)}\leq \mu(B)^{-1/2}.
$$
Define then $H^1(E)$ by the same procedure as for $H^1(\Gamma)$. \par
Our result is:
\begin{theo} \label{rieszfaible}
Assume that $\Gamma$ satisfies \eqref{doubling} and $(P_1)$. Then $d(I-P)^{-1/2}$ maps continuously $H^1(\Gamma)$ into $H^1(E)$.
\end{theo}
The proof goes through a duality argument. Let us introduce the $BMO(E)$ space. A function $\Phi$ on $E$ belongs to $BMO(E)$ if, and only if, $\Phi$ is antisymmetric and 
$$
\left\Vert \Phi\right\Vert_{BMO(E)}:=\left(\sup_{B\subset E} \frac 1{\mu(B)} \sum_{(x,y)\in B} \left\vert \Phi(x,y)-\Phi_B\right\vert^2d\mu_{xy}\right)^{1/2}<+\infty,
$$
where the supremum is taken over all balls $B\subset E$ and, as usual,
$$
\Phi_B:=\frac 1{\mu(B)} \sum_{(x,y)\in B} \Phi(x,y)\mu_{xy}.
$$
Define also $CMO(E)$ as the closure in $BMO(E)$ of the space of antisymmetric functions on $E$ with bounded support. Since $E$ is a space of homogeneous type, one has (\cite{coifmanweiss}):
\begin{theo} \label{dualhardy}
\begin{itemize}
\item[$1.$]
The dual of $H^1(E)$ is $BMO(E)$.
\item[$2.$]
The dual of $CMO(E)$ is $H^1(E)$.
\end{itemize}
\end{theo}
As in the proof of Theorem \ref{Riesz}, Theorem \ref{rieszfaible} will be a consequence of:
\begin{prop} \label{propatom}
Assume \eqref{doubling} and $(P_1)$. Then there exists $C>0$ such that, for all atom $a\in H^1(\Gamma)$,
$$
\left\Vert d(I-P)^{-1/2}a\right\Vert_{H^1(E)}\leq C.
$$
\end{prop}
{\bf Proof of Proposition \ref{propatom}: } we argue similarly to the proof of \cite{asterisque}, Chapter 4, Lemma 11 (see also Theorem 1 in \cite{mariasruss}), and will therefore be very sketchy. Let $a$ be an atom in $H^1(\Gamma)$ supported in a ball $B$. By assertion $2$ in Theorem \ref{dualhardy}, it is enough to prove that, for all antisymmetric function $\Phi$ on $E$ with bounded support,
\begin{equation} \label{testphi}
\left\vert \sum_{(x,y)\in E} d(I-P)^{-1/2}a(x,y)\Phi(x,y)\mu_{xy}\right\vert \leq C\left\Vert \Phi\right\Vert_{BMO(E)}.
\end{equation}
Since $d(I-P)^{-1/2}a\in L^1(E)$ and
$$
\sum_{(x,y)\in E} d(I-P)^{-1/2}a(x,y)\mu_{xy}=0,
$$
one has
\begin{equation} \label{identity}
\sum_{(x,y)\in E} d(I-P)^{-1/2}a(x,y)\Phi(x,y)\mu_{xy}=\sum_{(x,y)\in E} d(I-P)^{-1/2}a(x,y)\left(\Phi(x,y)-\Phi_{2B}\right)\mu_{xy},
\end{equation}
and \eqref{testphi} is derived from \eqref{identity} as in the proof of Lemma 11 in Chapter 4 of \cite{asterisque}. \hfill\fin\par

\medskip

\noindent Here is another result about the boundedness of Riesz transforms on Hardy spaces. A function $u:\Gamma\rightarrow \R$ is said to be harmonic on $\Gamma$ if and only if $(I-P)u(x)=0$ for all $x\in \Gamma$. Then:
\begin{theo} \label{rieszscalar}
Let $u:\Gamma\rightarrow \R$ be a harmonic function on $\Gamma$. Assume that there exist $x_0\in \Gamma$ and $C>0$ such that, for all $x\in \Gamma$,
$$
\left\vert u(x)\right\vert\leq C(1+d(x_0,x)).
$$
Define, for all functions $f$ on $\Gamma$ and all $x\in \Gamma$,
$$
R_u(f)(x)=\sum_{y\in \Gamma} d(I-P)^{-1/2}f(x,y)du(x,y)\mu_{xy}.
$$
Then $R_u$ is $H^1(\Gamma)$ bounded.
\end{theo}
Theorem \ref{rieszscalar} is a discrete counterpart of Theorem 1 in \cite{mariasruss} and the proof goes through a duality argument, as in the proof of Theorem 1 in \cite{mariasruss}. Indeed, the $H^1(\Gamma)-L^1(E)$ boundedness of $f\mapsto d(I-P)^{-1/2}f$ yields that $R_u$ is $H^1(\Gamma)-L^1(\Gamma)$ bounded. Then, if $f\in H^1(\Gamma)$, one checks that
\begin{equation} \label{cancel}
\sum_{x\in \Gamma} R_uf(x)m(x)=0.
\end{equation}
Indeed,
$$
\begin{array}{lll}
\displaystyle \sum_{x\in \Gamma} R_uf(x)m(x) &= & \displaystyle \sum_{x\in \Gamma}Êm(x)\sum_{y\sim x} d(I-P)^{-1/2}f(x,y)du(x,y)p(x,y)\\
& = & \displaystyle \sum_{y\in \Gamma}Ê\left(\sum_{x\in \Gamma} d(I-P)^{-1/2}f(x,y)du(x,y)p(x,y) m(x)\right)\\
& =& \displaystyle \sum_{x,y} d(I-P)^{-1/2}f(x,y)du(x,y) \mu_{xy}\\
& = & \displaystyle \langle d(I-P)^{-1/2}f,du\rangle_{L^2(E)}\\
& = & \displaystyle \langle (I-P)^{-1/2}f,\delta du\rangle_{L^2(E)}\\
& = & 0,
\end{array}
$$
since $\delta du=0$. Then, using \eqref{cancel}, one proves, arguing as in \cite{mariasruss}, that, if $a$ is an atom in $H^1(\Gamma)$, then, for all functions $\varphi$ with bounded support on $\Gamma$,
$$
\left\vert \sum_{x\in \Gamma}ÊR_uf(x)\varphi(x)m(x)\right\vert \lesssim \left\Vert \varphi\right\Vert_{BMO(\Gamma)}.
$$
The fact that $H^1(\Gamma)$ is the dual space of $CMO(\Gamma)$ then shows that
$$
\left\Vert R_ua\right\Vert_{H^1(\Gamma)}\leq C,
$$
and one concludes using the atomic decomposition for functions in $H^1(\Gamma)$. \hfill\fin\par

\medskip

\noindent Let us make a few comments on Theorems \ref{rieszfaible} and \ref{rieszscalar}. The conclusion of Theorem \ref{rieszfaible} says that, if $f\in H^1(\Gamma)$, then $d(I-P)^{-1/2}f$ has an atomic decomposition of the form
$$
d(I-P)^{-1/2}f=\sum_{k\in \N} \lambda_kA_k
$$
where $\sum_{k} \left\vert \lambda_k\right\vert\leq C\left\Vert f\right\Vert_{H^1(\Gamma)}$ and the $A_k$'s are atoms in $H^1(E)$. However, one does not claim that each $A_k$ is equal to $da_k$ where $a_k$ is an atom in $\dot{S}^{1,1}(\Gamma)$. In this sense, the conclusion of Theorem \ref{rieszfaible} is weaker than the one of Theorem \ref{Riesz}. On the other hand, assumption \eqref{reversevol} is not required in Theorem \ref{rieszfaible}. Finally, Theorem \ref{rieszscalar} says that a scalar version of the Riesz transform is $H^1(\Gamma)$-bounded and does not require assumption \eqref{reversevol} either. \par

{\bf Acknowledgements: } the authors would like to thank G. Dafni and E. M. Ouhabaz for useful remarks on this manuscript.
\bibliographystyle{alpha}               
  \bibliography{HSA.bib}   

\begin{thebibliography}{DMRT10}

\bibitem[AC05]{auscher2005riesz}
P.~Auscher and T.~Coulhon.
\newblock {Riesz transform on manifols and Poincar{\'e} inequalities}.
\newblock {\em Ann. Sc. Norm. Super. Pisa Cl. Sci. (5)}, 4(3):531--555, 2005.

\bibitem[AMR08]{amr}
P.~Auscher, A.~McIntosh, and E.~Russ.
\newblock Hardy spaces of differential forms on {R}iemannian manifolds.
\newblock {\em J. Geom. Anal.}, 18(1):192--248, 2008.

\bibitem[ART05]{art}
P.~Auscher, E.~Russ, and P.~Tchamitchian.
\newblock {Hardy Sobolev spaces on strongly Lipschitz domains of ${\mathbb
  R}^n$}.
\newblock {\em J. Funct. Anal.}, 218(1):54--109, 2005.

\bibitem[AT98]{asterisque}
P.~Auscher and P.~Tchamitchian.
\newblock {Square root problem for divergence operators and related topics}.
\newblock {\em Asterisque}, 249, 1998.

\bibitem[BB10]{badber}
N.~Badr and F.~Bernicot.
\newblock Abstract {H}ardy-{S}obolev spaces and interpolation.
\newblock {\em J. Funct. Anal.}, 259(5):1169--1208, 2010.

\bibitem[BD10]{badr2010atomic}
N.~Badr and G.~Dafni.
\newblock An atomic decomposition of the {H}aj\l asz {S}obolev space {$M^1_1$}
  on manifolds.
\newblock {\em J. Funct. Anal.}, 259(6):1380--1420, 2010.

\bibitem[BD11]{maxhs}
N.~Badr and G.~Dafni.
\newblock Maximal characterization of hardy-sobolev spaces on manifolds.
\newblock In {\em Concentration, Functional Inequalities and Isoperimetry},
  volume 545 of {\em Contemp. Math.}, pages 13--21. Amer. Math. Soc., 2011.

\bibitem[BR09]{BR09}
N.~Badr and E.~Russ.
\newblock Interpolation of {S}obolev spaces, {L}ittlewood-{P}aley inequalities
  and {R}iesz transforms on graphs.
\newblock {\em Publ. Mat.}, 53(2):273--328, 2009.

\bibitem[Cal72]{calderon}
A.~P. Calder{\'o}n.
\newblock Estimates for singular integral operators in terms of maximal
  functions.
\newblock {\em Studia Math.}, 44:563--582, 1972.
\newblock Collection of articles honoring the completion by Antoni Zygmund of
  50 years of scientific activity, VI.

\bibitem[CG98]{cougri}
T.~Coulhon and A.~Grigoryan.
\newblock {Random walks on graphs with regular volume growth}.
\newblock {\em Geom. Funct. Anal.}, 8(4):656--701, 1998.

\bibitem[CGZ05]{cgz}
T.~Coulhon, A.~Grigor'yan, and F.~Zucca.
\newblock The discrete integral maximum principle and its applications.
\newblock {\em Tohoku Math. J. (2)}, 57(4):559--587, 2005.

\bibitem[CW77]{coifmanweiss}
R.~R. Coifman and G.~Weiss.
\newblock Extensions of {H}ardy spaces and their use in analysis.
\newblock {\em Bull. Amer. Math. Soc.}, 83(4):569--645, 1977.

\bibitem[DMRT10]{dmrt}
R.~Duran, M.-A. Muschietti, E.~Russ, and P.~Tchamitchian.
\newblock Divergence operator and {P}oincar\'e inequalities on arbitrary
  bounded domains.
\newblock {\em Complex Var. Elliptic Equ.}, 55(8-10):795--816, 2010.

\bibitem[FS72]{feffstein}
C.~Fefferman and E.~M. Stein.
\newblock {$H^{p}$} spaces of several variables.
\newblock {\em Acta Math.}, 129(3-4):137--193, 1972.

\bibitem[Haj96]{haj}
P.~Haj{\l}asz.
\newblock Sobolev spaces on an arbitrary metric space.
\newblock {\em Potential Anal.}, 5(4):403--415, 1996.

\bibitem[Haj03a]{hajlasznew}
P.~Haj{\l}asz.
\newblock A new characterization of the {S}obolev space.
\newblock {\em Studia Math.}, 159(2):263--275, 2003.
\newblock Dedicated to Professor Aleksander Pe{\l}czy{\'n}ski on the occasion
  of his 70th birthday (Polish).

\bibitem[Haj03b]{hajlasz}
P.~Haj{\l}asz.
\newblock Sobolev spaces on metric-measure spaces.
\newblock In {\em Heat kernels and analysis on manifolds, graphs, and metric
  spaces ({P}aris, 2002)}, volume 338 of {\em Contemp. Math.}, pages 173--218.
  Amer. Math. Soc., Providence, RI, 2003.

\bibitem[HK98]{quasicont}
P.~Haj{\l}asz and J.~Kinnunen.
\newblock H\"older quasicontinuity of {S}obolev functions on metric spaces.
\newblock {\em Rev. Mat. Iberoamericana}, 14(3):601--622, 1998.

\bibitem[HK00]{sobmetpoinc}
P.~Haj{\l}asz and P.~Koskela.
\newblock Sobolev met {P}oincar\'e.
\newblock {\em Mem. Amer. Math. Soc.}, 145(688):x+101, 2000.

\bibitem[KS08]{kosksaks}
P.~Koskela and E.~Saksman.
\newblock Pointwise characterizations of {H}ardy-{S}obolev functions.
\newblock {\em Math. Res. Lett.}, 15(4):727--744, 2008.

\bibitem[KT07]{kin2007}
J.~Kinnunen and H.~Tuominen.
\newblock Pointwise behaviour of {$M^{1,1}$} {S}obolev functions.
\newblock {\em Math. Z.}, 257(3):613--630, 2007.

\bibitem[KZ08]{kz}
S.~Keith and X.~Zhong.
\newblock The poincar\'e inequality is an open ended condition.
\newblock {\em Ann. Math.}, 167:575--599, 2008.

\bibitem[Mar01]{martell01}
J.M. Martell.
\newblock {\em {Desigualdades con pesos en el Analisis de Fourier: de los
  espacios de tipo homogeneo a las medidas no doblantes}}.
\newblock PhD thesis, Ph. D. Thesis, Universidad Autonoma de Madrid, 2001.

\bibitem[Miy90]{miyachi}
A.~Miyachi.
\newblock Hardy-{S}obolev spaces and maximal functions.
\newblock {\em J. Math. Soc. Japan}, 42(1):73--90, 1990.

\bibitem[MR03]{mariasruss}
M.~Marias and E.~Russ.
\newblock {$H^1$-boundedness of Riesz transforms and imaginary powers of the
  Laplacian on Riemannian manifolds}.
\newblock {\em Ark. Mat.}, 41(1):115--132, 2003.

\bibitem[Rus00]{scand}
E.~Russ.
\newblock {Riesz transforms on graphs for $1\leq p\leq 2$}.
\newblock {\em Math. Scand.}, 87:133--160, 2000.

\bibitem[Rus01]{pota}
E.~Russ.
\newblock {$H^1-L^1$ Boundedness of Riesz Transforms on Riemannian Manifolds
  and on Graphs}.
\newblock {\em Pot. Anal.}, 14(3):301--330, 2001.

\bibitem[Ste93]{E.M.Stein93}
E.~M. Stein.
\newblock {\em Harmonic analysis: real-variable methods, orthogonality, and
  oscillatory integrals}, volume~43 of {\em Princeton Mathematical Series}.
\newblock Princeton University Press, Princeton, NJ, 1993.

\bibitem[Str90]{stri}
R.~S. Strichartz.
\newblock {$H^p$} {S}obolev spaces.
\newblock {\em Colloq. Math.}, 60/61(1):129--139, 1990.

\bibitem[SW71]{steinweiss}
E.~M. Stein and G.~Weiss.
\newblock {\em Introduction to {F}ourier analysis on {E}uclidean spaces}.
\newblock Princeton University Press, Princeton, N.J., 1971.
\newblock Princeton Mathematical Series, No. 32.

\end{thebibliography}

\end{document}